\documentstyle{article}

\newtheorem{theorem}{Theorem}[section]
\newtheorem{remark}[theorem]{Remark}
\newtheorem {lemma}[theorem]{Lemma}
\newtheorem {corollary}[theorem]{Corollary}
\newtheorem {proposition}[theorem]{Proposition}
\newtheorem{definition}[theorem]{Definition}
\newtheorem{example}[theorem]{Example}

\newtheorem{assumption}[theorem]{Assumption}

\def\cC{\mathcal{C}}
\def\cF{\mathcal{F}}
\def\cE{\mathcal{E}}
\def\cB{\mathcal{B}}

\def\cS{\mathcal{S}}

\newdimen\AAdi%
\newbox\AAbo%
\def\AArm{\fam0 }
\def\AAk#1#2{\setbox\AAbo=\hbox{#2}\AAdi=\wd\AAbo\kern#1\AAdi{}}%
\def\AAr#1#2#3{\setbox\AAbo=\hbox{#2}\AAdi=\ht\AAbo\raise#1\AAdi\hbox{#3}}%

\def\BBn{{\AArm I\!N}}%
\def\BBr{{\AArm I\!R}}%
\def\BBz{{\AArm Z\!\!Z}}%

\begin{document}

\title{Convergence of Point Processes with Weakly Dependent Points}

\author{Raluca M. Balan\footnote{Corresponding author. Department of Mathematics and Statistics, University of Ottawa,
585 King Edward Avenue, Ottawa, ON, K1N 6N5, Canada. E-mail
address: rbalan@uottawa.ca} \ \footnote{Supported by a grant from
the Natural Sciences and Engineering Research Council of Canada.}
\ and Sana Louhichi\footnote{Laboratoire de Probabilit\'es,
Statistique et mod\'elisation, Universit\'e de Paris-Sud, B\^at.
425, F-91405 Orsay Cedex, France. E-mail address:
Sana.Louhichi@math.u-psud.fr}}

\date{May 19, 2008}

\maketitle

\begin{abstract}

\noindent For each $n \geq 1$, let $\{X_{j,n}\}_{1 \leq j \leq n}$
be a sequence of strictly stationary random variables. In this
article, we give some asymptotic weak dependence conditions for
the convergence in distribution of the point process
$N_n=\sum_{j=1}^{n}\delta_{X_{j,n}}$ to an infinitely divisible
point process. From the point process convergence, we obtain the
convergence in distribution of the partial sum sequence
$S_n=\sum_{j=1}^{n}X_{j,n}$ to an infinitely divisible random
variable, whose L\'{e}vy measure is related to the canonical
measure of the limiting point process. As examples, we discuss the
case of triangular arrays which possess known (row-wise)
dependence structures, like the strong mixing property, the
association, or the dependence structure of a stochastic
volatility model.

\end{abstract}

{\em Keywords:} weak limit theorem; point process; infinitely
divisible law; strong mixing; association


\section{Introduction} 

In this article we examine the convergence in distribution of the
sequence $\{N_n=\sum_{j=1}^{n}\delta_{X_{j,n}},n \geq 1\}$ of
point processes, whose points $(X_{j,n})_{1 \leq j \leq n, n \geq
1}$ form a triangular array of strictly stationary weakly
dependent random variables, with values in a locally compact
Polish space $E$. As it is well-known in the literature, by a mere
invocation of the continuous mapping theorem, the convergence of
the point processes $(N_n)_{n \geq 1}$ becomes a rich source for
numerous other limit theorems, which describe the asymptotic
behavior of various functions of $N_n$, provided these functions
are continuous with respect to the vague topology (in the space of
Radon measures in which $N_n$ lives). This turns out to be a very
useful approach, provided that one has a handle on the limiting
point process $N$, which usually comes  from its ``cluster
representation'', namely a representation of the form $N
\stackrel{d}{=}\sum_{i,j \geq 1}\delta_{T_{ij}}$ for some
carefully chosen random points $T_{ij}$. In principle, one can
obtain via this route the convergence of the partial sum sequence
$\{S_n=\sum_{j=1}^{n}X_{j,n}, n \geq 1\}$ to the sum $X:=\sum_{i,j
\geq 1}T_{ij}$ of the points. (However, as it is usually the case
in mathematics, this works only ``in principle'', meaning that the
details are not to be ignored.)

On the other hand, a classical result in probability theory says
that the class of all limiting distributions for the partial sum
sequence $(S_n)_{n \geq 1}$ associated with a triangular array of
independent random variables coincides with the class of all
infinitely divisible distributions (see e.g. Theorem 4.2,
\cite{gnedenko-kolmogorov54}). This result has been extended
recently in \cite{dedecker-louhichi05a} and
\cite{dedecker-louhichi05b} to some similar results for arrays of
weakly dependent random variables with finite variances. One of
the goals of the present article is to investigate if such a limit
theorem can be obtained via the more powerful approach of point
process convergence, which does not require any moment
assumptions.

Our work is a continuation of the line of research initiated by
Davis and Hsing in their magistral article \cite{davis-hsing95},
in which they consider an array of random variables with values in
$\BBr \verb2\2 \{0\}$, of the form $X_{j,n}=X_j/a_n$, where
$(X_j)_{j \geq 1}$ is a strictly stationary sequence with heavy
tails, and $a_n$ is the $(1-1/n)$-quantile of $X_1$. In this case,
there is no surprise that the limiting distribution of the
sequence $(S_n)_{n \geq 1}$ coincides with the stable law.

The main asymptotic dependence structure in our array (called
(AD-1)) is inherited from condition ${\cal A}(\{a_n\})$ of
\cite{davis-hsing95} (see also \cite{davis-mikosch98}), but unlike
these authors, we do not require that $X_{1,n}$ lie in the domain
of attraction of the stable law. This relaxation will allow us
later to obtain a general (possibly non-stable) limit distribution
for the partial sum sequence $(S_n)_{n \geq 1}$. The asymptotic
dependence structure (AD-1) is based on a simple technique, which
requires that we ``split'' each row  of the array into $k_n$
blocks of length $r_n$, and then we ask the newly produced
block-vectors $Y_{i,n}=(X_{(i-1)r_n+1,n}, \ldots, X_{ir_n,n}), 1
\leq i \leq k_n$ to behave asymptotically as their independent
copies $\tilde Y_{i,n}=(\tilde X_{(i-1)r_n+1,n}, \ldots, \tilde
X_{ir_n,n}), 1 \leq i \leq k_n$. Note that this procedure does not
impose any restrictions on the dependence structure within the
blocks, it only specifies (asymptotically) the dependence
structure {\em among} the blocks. (Since each row of the array has
length $n$, as a by-product, this procedure necessarily yields a
remainder number $n-r_n k_n$ of terms, which will be taken care of
by an asymptotic negligibility condition called (AN).) The origins
of this technique can be traced back to Jakubowski's thorough
investigation of the minimal asymptotic dependence conditions in
the stable limit theorems (see \cite{jakubowski93},
\cite{jakubowski97}). However, in Jakubowski's condition (B) the
number of blocks is assumed to be $2$, whereas in our condition
(AD-1), as well as in condition ${\cal A}(\{a_n\})$, the number
$k_n$ of blocks explodes to infinity.

A fundamental result of \cite{davis-hsing95} states that under
${\cal A}(\{a_n\})$, {\em if} the sequence of point processes
$(N_n)_{n \geq 1}$ converges, then its limit $N$ admits a very
nice cluster representation of the form $N=\sum_{i,j \geq
1}\delta_{P_i Q_{ij}}$, with independent components $(P_i)_{i \geq
1}$ and $(Q_{ij})_{i \geq 1,j \geq 1}$ (see Theorem 2.3,
\cite{davis-hsing95}). The key word in this statement is ``if''.
In the present article, we complement Theorem 2.3,
\cite{davis-hsing95} by supplying a new asymptotic dependence
condition (called (AD-2)) which along with conditions (AD-1) and
(AN), ensure that the convergence of $(N_n)_{n \geq 1}$ {\em does}
happen, for an arbitrary triangular array of random variables (not
necessarily of the form $X_{j,n}=X_j/a_n$, with $a_n \sim
n^{1/\alpha}L(n)$ for some $ \alpha \in (0,2)$ and a slowly
varying function $L$). Under this new condition, we are able to
find a new formulation for Kallenberg's necessary and sufficient
condition for the convergence of $(N_n)_{n \geq 1}$, in terms of
the incremental differences between the Laplace functionals of the
processes $N_{m,n}=\sum_{j=1}^{m}\delta_{X_{j,n}}, m \leq n$.

The new condition (AD-2) is an ``anti-clustering'' condition,
which does not allow the jump times of the partial sum process
$S_n(t)=\sum_{j=1}^{[nt]}X_{j,n}, t \in [0,1]$, whose size exceed
in modulus an arbitrary fixed value $\eta>0$ (and are located at a
minimum distance of $m/n$ of each other), to get condensed in a
``small'' interval of time of length $r_n/n \sim k_n^{-1}$. This
happens with a probability which is asymptotically $1$, when $n$
gets large and $m$ either stabilizes around a finite value $m_0$,
or gets large as well. From the mathematical point of view,
condition (AD-2) treats the dependence structure {\em within} the
blocks of length $r_n$, which was left open by condition (AD-1).
Condition (AD-2) is automatically satisfied when the rows of the
array are $m$-depedent. The asymptotic negligibility condition
(AN) forces the rate of the convergence in probability to $0$ of
$X_{1,n}$ to be at most $n^{-1}$.

Of course, there are many instances in the literature in which the
sequence $(N_n)_{n \geq 1}$ converges.
The most notable example is probably the case of the moving
average sequences
$X_{i,n}=a_n^{-1}\sum_{j=0}^{\infty}C_{i,j}Z_{i-j}, i \geq 1$:
the classical Theorem 2.4.(i), \cite{davis-resnick85} treats the
case of constant coefficients $C_{i,j}=c_{j}$, whereas the recent
Theorem 3.1, \cite{kulik06} allows for random coefficients
$C_{i,j}$. Another important example is given by Theorem 3.1,
\cite{resnick-samorodnitsky04}, in which
$X_{j,n}=n^{-1/\alpha}X_j$ and $(X_j)_{j \geq 1}$ is a symmetric
$\alpha$-stable process.



With the convergence of the sequence $(N_n)_{n \geq 1}$ in hand,
we can prove a general (non-stable) limit theorem for the partial
sum sequence $(S_n)_{n \geq 1}$, and the hypothesis of this new
theorem are indeed verified by the moving average sequences (even
with random coefficients).
The infinitely divisible law that we obtain as the limit of
$(S_n)_{n \geq 1}$ must have a L\'evy measure $\rho$ which
satisfies the condition $\int_0^1 x \rho(dx)<\infty$. This is a
limitation which has to do with the method that we use, based on
the Ferguson-Klass \cite{ferguson-klass72} representation of an
infinitely divisible law. It remains an open problem to see how
one could recover a general infinitely divisible law as the limit
of $(S_n)_{n \geq 1}$, using point process techniques.

There is a large amount of literature dedicated to limit theorems
for the partial sum sequence associated to a triangular array,
based on point process techniques. For a comprehensive account on
this subject in the independent case, we refer the reader to the
expository article \cite{resnick86}. In the case of arrays which
possess a row-wise dependence structure, the first systematic
application of point process techniques for obtaining limit
theorems for the sequence $(S_n)_{n \geq 1}$, has been developed
in \cite{jakubowski-kobus89}. The article
\cite{jakubowski-kobus89} identifies the necessary conditions for
the general applicability of point process techniques, including
cases which are not covered in the present article (e.g. the case
of the $\alpha$-stable limit distribution, with $\alpha \in
[1,2)$), and contains the first general limit theorem for sums of
$m$-dependent random variables with heavy tails. Without aiming at
exhausting the entire list of contributions to this area, we
should also mention the article \cite{kobus95}, which includes the
necessary and sufficient conditions for the convergence in
distribution of sums of $m$-dependent random variables, to a
generalized Poisson distribution.

The present article is organized as follows. In Section 2, we
introduce the asymptotic dependence conditions and we prove the
main theorem which gives the convergence of the sequence $(N_n)_{n
\geq 1}$ of point processes. Section 3 is dedicated to the
convergence of the partial sum sequence $(S_n)_{n \geq 1}$. In
Section 4, we give a direct consequence of the main theorem, which
can be viewed as a complement of Theorem 2.3,
\cite{davis-hsing95}, specifying some conditions which guarantee
that the limit $N$ of the sequence $(N_n)_{n \geq 1}$ exists (and
admits a ``product-type'' cluster representation). Section 5 is
dedicated to the analysis of condition (AD-1) in the case of an
array whose row-wise dependence structure is given by the strong
mixing property, the association, or is that of a stochastic
volatility sequence. Appendix A gives a necessary and sufficient
condition for a product-type cluster representation of a Poisson
process. Appendix B gives a technical construction needed in the
strongly mixing case and in the case of a stochastic volatility
sequence.

\section{Weak convergence of point processes}
\label{weak-convergence}

We begin by introducing the point process background. Our main
references are \cite{kallenberg83}, \cite{resnick87} and
\cite{resnick92}.

Let $E$  be a locally compact Polish space, $\cE$ its Borel
$\sigma$-algebra and $\cB$ the class of bounded Borel sets in
$\cE$. A measure $\mu$ on $E$ is called Radon if $\mu(B)<\infty$
for all $B \in \cB$. If $E=\BBr \verb2\2 \{0\}$ or $E=(0,\infty)$,
the class $\cB$ contains the Borel sets in $E$ which are bounded
away from $0$ and $\pm \infty$, respectively from $0$ and
$\infty$.

 Let
$M_p(E)$ be the class of all Radon measures on $E$ such that
$\mu(B) \in \BBz_{+}=\{0,1,2, \ldots\}$ for all $B \in \cB$. The
space $M_p(E)$ is endowed with the topology of vague convergence.
The corresponding Borel $\sigma$-field is denoted by ${\cal
M}_p(E)$. (Recall that a sequence $(\mu_n)_{n \geq 1} \subset
M_{p}(E)$ converges vaguely to $\mu$ if $\mu_n(B) \to \mu(B)$ for
any $B \in \cB$ with $\mu(\partial B)=0$.) For each
$B_1,\ldots,B_k \in \cE$, we define  $\pi_{B_1,\ldots, B_k}:
M_{p}(E) \to \BBz_{+}^k$  by $\pi_{B_1,\ldots,B_k}(\mu)=(\mu(B_1),
\ldots, \mu(B_k))$. We denote by $\delta_x$ the Dirac measure at
$x \in E$, and by $o$ the null measure in $M_p(E)$. For any $\mu
\in M_p(E)$ and for any measurable non-negative function $f$ on
$E$, we let $\mu(f)=\int_E f(x)\mu(dx)$.

A point process $N$ is an $M_p(E)$-valued random variable, defined
on a probability space $(\Omega,\cF,P)$. Its Laplace functional is
defined by $L_{N}(f)=E(e^{-N(f)})$, for any measurable
non-negative function $f$ on $E$. If $N_1$ and $N_2$ are two point
processes on the same probability space, we use the notation $N_1
\stackrel{d}{=} N_2$ if $P \circ N_1^{-1}=P \circ N_2^{-1}$; this
is equivalent to the fact that $L_{N_1}(f)=L_{N_2}(f)$, for any
measurable non-negative function $f$ on $E$.

If $N,(N_n)_{n \geq 1}$ are point processes, we say that $(N_n)_{n
\geq 1}$ converges in distribution to $N$ (and we write $N_n
\stackrel{d}{\rightarrow} N$), if $\{P \circ N_n^{-1}\}_{n \geq
1}$ converges weakly to $P \circ N^{-1}$.
By the continuous mapping theorem, if $N_n
\stackrel{d}{\rightarrow} N$, then $\{h(N_n)\}_{n \geq 1}$
converges in distribution to $h(N)$, for every continuous function
$h: M_p(E) \to \BBr$. By Theorem 4.2, \cite{kallenberg83},
$(N_{n})_{n \geq 1}$ converges in distribution to $N$ if and only
if $L_{N_n}(f) \to L_{N}(f), \forall f \in C_{K}^{+}(E)$, where
$C_{K}^{+}(E)$ denotes the class of all continuous non-negative
functions $f$ on $E$, with compact support.

A point process $N$ is said to be infinitely divisible if for
every $n \geq 1$, there exist some i.i.d. point processes $N_1,
\ldots, N_n$ such that $N \stackrel{d}{=}N_1+ \ldots + N_n$. By
Theorem 6.1, \cite{kallenberg83}, if $N$ is an infinitely
divisible process, then there exists a unique measure $\lambda$ on
$M_{p}(E) \verb2\2 \{o\}$ (called the canonical measure of $N$)
such that
\begin{equation}
\label{cond-lambda-1}\int_{M_{p}(E) \verb2\2 \{o\}}
(1-e^{-\mu(B)})\lambda(d\mu)<\infty, \quad \forall B \in \cB \
\mbox{and}
\end{equation}
\begin{equation}
\label{log-Laplace-N} -\log L_{N}(f)=\int_{M_{p}(E) \verb2\2
\{o\}}(1-e^{-\mu(f)})\lambda(d\mu), \quad \forall f \in
C_{K}^{+}(E).
\end{equation}

We begin to introduce our framework. For each $n \geq 1$, let
$(X_{j,n})_{1\leq j\leq n}$ be a strictly stationary sequence of
$E$-valued random variables, defined on a probability space
$(\Omega,\cF, P)$.

We introduce our first asymptotic dependence condition. A similar
condition was considered in \cite{davis-hsing95}.

\begin{definition}
We say that the triangular array $(X_{j,n})_{1 \leq j \leq n,  n
\geq 1}$ satisfies {\bf condition (AD-1)} if there exists a
sequence $(r_n)_n \subset \BBz_+$ with $r_n \to \infty$ and
$k_n=[n/r_n] \to \infty$ as $n \to \infty$, such that:
\begin{equation}
\label{cond-AD-1} \lim_{n\rightarrow \infty}\left|E \left(
e^{-\sum_{j=1}^{n} f(X_{j,n})} \right)-
\left\{E\left(e^{-\sum_{j=1}^{r_n}f(X_{j,n})}\right)
\right\}^{k_n} \right|=0, \quad \forall f \in C_K^+(E).
\end{equation}
\end{definition}

In Section \ref{section-about-AD}, we will examine condition
(AD-1) in the case of arrays which possess a known dependence
structure on each row.

To see the intuitive meaning of condition (AD-1), let us consider
the point process $N_{m,n}=\sum_{j=1}^{m}\delta_{X_{j,n}}$, whose
Laplace functional is denoted by $L_{m,n}$, for each $m \leq n$.
By convention, we let $L_{0,n}=1$. We denote $N_n=N_{n,n}$. Note
that
$L_{m,n}(f)=E(e^{-N_{m,n}(f)})=E(e^{-\sum_{j=1}^{m}f(X_{j,n})})$.

For each $n \geq 1$, let $(\tilde N_{i,n})_{1 \leq i \leq k_n}$ be
a sequence of i.i.d. point processes with the same distribution as
$N_{r_{n},n}$, and let $\tilde N_n=\sum_{i=1}^{k_n}\tilde
N_{i,n}$. Then $L_{\tilde N_n}(f)=\{L_{r_n,n}(f)\}^{k_n}$ and
(\ref{cond-AD-1}) becomes: $|L_{N_n}(f) -L_{\tilde N_n}(f)| \to 0,
\forall f \in C_K^+(E)$.

\noindent This shows that under (AD-1), the asymptotic behavior of
the sequence $(N_n)_n$ is the same as that of $(\tilde N_n)_n$,
i.e. $(N_n)_n$ converges in distribution if and only if $(\tilde
N_n)_n$ does, and in this case, the limits are the same.

\vspace{2mm}

We now introduce an asymptotic negligibility condition, in
probability.

\begin{definition}
We say that $(X_{j,n})_{j \leq n, n \geq 1}$ satisfies {\bf
condition (AN)} if $$\limsup_{n \to \infty}nP(X_{1,n} \in
B)<\infty, \quad \forall B \in \cB.$$
\end{definition}

Under (AN), the triangular array $(\tilde N_{i,n})_{1 \leq i \leq
k_n, \ n \geq 1}$ becomes a ``null-array'', i.e. $P(\tilde
N_{1,n}(B) >0) \to 0$ for all $B \in \cB$. To see this, note that
$P(\tilde N_{1,n}(B)>0)=P(\bigcup_{j=1}^{r_n} \{X_{j,n} \in B\})
\leq (n/k_n)P(X_{1,n} \in B) \to 0$. By invoking Theorem 6.1,
\cite{kallenberg83}, we infer that the sequence $(\tilde N_n)_{n
\geq 1}$ (or equivalently, the sequence $(N_n)_{n \geq 1}$)
converges in distribution to {\em some} point process $N$ if and
only if there exists a measure $\lambda$ satisfying
(\ref{cond-lambda-1}) such that
\begin{equation}
\label{Kallenberg-condition} k_n(1-L_{r_n,n}(f)) \to \int_{M_p(E)
\verb2\2 \{o\}}(1-e^{-\mu(f)})\lambda(d\mu), \quad \forall f \in
C_{K}^{+}(E).
\end{equation}
In this case, $N$ is an infinitely divisible point process with
canonical measure $\lambda$, i.e. (\ref{log-Laplace-N}) holds. By
writing
$$L_{\tilde N_n}(f)=\{L_{r_n,n}(f)\}^{k_n}=\left\{1-\frac{k_n(1-L_{r_n,n}(f))}{k_n}
\right\}^{k_n}$$ and using the fact that $(1+x_n/n)^n \to e^x$ iff
$x_n \to x$, we see that condition (\ref{Kallenberg-condition})
requires that $L_{\tilde N_n}(f) \to L_{N}(f)$.

In conclusion, when dealing with triangular arrays which satisfy
(AD-1) and (AN), the {\em only} possible limit (if it exists) for
the sequence $(N_n)_{n \geq 1}$ of point processes is an
infinitely divisible point process.

As in \cite{dedecker-louhichi05a}, if $(r_n)_n$ is an arbitrary
sequence of positive integers with $r_n \to \infty$, we let
${\cS}={\cS}_{(r_n)_n}$ be the set of all positive integers $m$
such that: $\limsup_{n \to \infty}n
\sum_{j=m+1}^{r_n}E[f(X_{1,n})f(X_{j,n})]=0, \ \forall f \in
C_K^{+}(E)$. Let $m_0$ be the smallest integer in $\cS$. By
convention, we let $m_0=\infty$ if $\cS=\emptyset$. For an
arbitrary function $\phi$, we denote
$$\lim_{m \to m_0}\phi(m)=
\left\{
\begin{array}{ll}
\phi(m_0), & \mbox{if $m_0<\infty$} \\
\lim_{m \to \infty}\phi(m),& \mbox{if $m_0=\infty$}
\end{array}
\right.$$

We are now ready to introduce our second asymptotic dependence
condition.

\begin{definition}
We say that the triangular array $(X_{j,n})_{1 \leq j \leq n, n
\geq 1}$ satisfies {\bf condition (AD-2)} if there exists a
sequence $(r_n)_n \subset {\BBz}_{+}$ with $r_n \to \infty$ and
$m_0:=\inf {\cS}_{(r_n)_n}$, such that
$$\lim_{m \to m_0}\limsup_{n \to \infty}n\sum_{j=m+1}^{r_n}
E[f(X_{1,n})f(X_{j,n})]=0, \quad \forall f \in C_{K}^{+}(E).$$
\end{definition}

Specifying the row-wise dependence structure of the array does not
guarantee that condition (AD-2) is satisfied, but it may help to
understand its meaning.

\begin{example}
{\rm  ($m$-dependent random variables) Suppose that for every $n
\geq 1$,  the sequence $(X_{j,n})_{1 \leq j \leq n}$ is {\em
$m$-dependent}, i.e. $(X_{1,n}, \ldots, X_{j,n})$ and $(X_{j+r,n},
X_{j+r+1,n}, \ldots, X_{n,n})$ are independent, for all $j,r \leq
n$ with $j+r \leq n$ and $r \geq m$. Suppose that the array
$(X_{j,n})_{1 \leq j \leq n, n \geq 1}$ satisfies condition (AN).
Then for any sequence $(r_n)_n \subset {\BBz}_{+}$ with $r_n \to
\infty$ and $k_n:=[n/r_n] \to \infty$, ${\cS}_{(r_n)_n}=\{l \in
{\BBz}_{+}; l \geq m\}$, since for any $l \geq m $ and for any $f
\in C_{K}^{+}(E)$,
$$n\sum_{j=l+1}^{r_n}E [f(X_{1,n})f(X_{j,n})] \leq n r_n \{E[
f(X_{1,n})]\}^2 \leq k_n r_n^2 \|f\|_{\infty}^2 P(X_{1,n}\in K)^2
\leq \frac{C}{k_n} \|f\|_{\infty}^2  \to 0.$$

\noindent (Here $K$ is the compact support of $f$.) Therefore
$m_0:=\inf {\cS}_{(r_n)_n}=m$ and condition (AD-2) is satisfied.
In particular, if the sequence $(X_{j,n})_{1 \leq j \leq n}$ is
$1$-dependent (or i.i.d.), then ${\cS}_{(r_n)_n}=\BBz_+$ and
$m_0=1$.}
\end{example}

\begin{remark}
{\rm The following slightly stronger form of condition (AD-2) has
a clearer intuitive meaning. We say that the triangular array
$(X_{j,n})_{1 \leq j \leq n, n \geq 1}$ satisfies {\bf condition
(AD-2')} if there exists a sequence $(r_n)_n \subset \BBz_{+}$
with $r_n \to \infty$ and $m_0:=\inf {\cS}_{(r_n)_n}$, such that
$$\lim_{m \to m_0}\limsup_{n \to \infty}n \sum_{j=m+1}^{r_n}
P(X_{1,n} \in B, X_{j,n} \in B)=0, \quad \forall B \in \cB.$$

\noindent Note that, due to the stationarity of the array, we
have:
\begin{equation}
\label{AD-2-stationarity} P\left(\bigcup_{i=1}^{r_n-m}
\bigcup_{k=m+i}^{r_n} \{X_{i,n} \in B, X_{k,n} \in B \}\right)
\leq
r_n \sum_{j=m+1}^{r_n} P(X_{1,n} \in B, X_{j,n} \in B).
\end{equation}

\noindent Therefore, if condition (AD-2') holds, and we let
$k_n:=[n/r_n]$, then
$$\lim_{m \to m_0} \limsup_{n \to \infty}k_n P(\exists \ i<k \leq r_n \ \mbox{with} \
k-i \geq m \ \mbox{such that} \
 X_{i,n} \in B, X_{k,n} \in B)= 0.$$

\noindent In particular, if condition (AD-2') holds with $m_0=1$,
then
$$k_n P(N_{r_n,n}(B) >1 )=k_n P( \exists \ i<k \leq r_n
\ \mbox{such that} \ X_{i,n} \in B, X_{k,n} \in B) \to 0.$$

\noindent For each $B \in \cB$ and for each $t \in [0,1]$, define
$M_n^{B}([0,t])=N_{[nt],n}(B)$. Condition (AD-2') with $m_0=1$
forces
$$\lim_{n \to \infty}\frac{n}{r_n}P\left(M_n^B \left(\left[0, \frac{r_n}{n}\right]\right)>1
\right) =0.$$ Intuitively, if $k_n \to \infty$, we can view this
as an ``asymptotic orderly'' property of the sequence
$(M_n^B)_{n}$. (According to p. 30, \cite{daley-verajones03}, a
point process $N$ is called {\em orderly} if $\lim_{t \to
0}t^{-1}P(N([0,t])>1)=0$.) }
\end{remark}

The following theorem gives a necessary and sufficient condition
for the convergence in distribution of the sequence $(N_n)_n$. As
mentioned earlier, the limit process must be an infinitely
divisible point process.

As it was pointed out by an anonymous referee, our approach to
identify the limit in the theorem below, is closely related to the
method used in the proof of Theorem 3.1 of \cite{jakubowski97}.

\begin{theorem}
\label{main-theorem} For each $n \geq 1$, let $(X_{j,n})_{1 \leq j
\leq n}$ be a strictly stationary sequence of $E$-valued random
variables. Suppose that the triangular array  $(X_{j,n})_{1 \leq j
\leq n, n \geq 1}$ satisfies condition (AN), as well as conditions
(AD-1) and (AD-2) (with the same sequence $(r_n)_n)$). Denote
$m_0:=\inf {\cS}_{(r_n)_n}$.

Then the sequence $(N_n)_{n \geq 1}$ converges in distribution to
some point process $N$ if and only if there exists a measure
$\lambda$ on $M_p(E)\setminus\{o\}$ which satisfies
(\ref{cond-lambda-1}), such that
\begin{equation}
\label{cond-L} \lim_{m\rightarrow m_0}\limsup_{n\rightarrow
\infty} \left|n(L_{m-1,n}(f)-L_{m,n}(f))-\int_{M_p(E) \verb2\2
\{o\}} (1-e^{-\mu(f)})\lambda(d\mu)\right|=0, \quad \forall f\in
C_K^+(E).
\end{equation}
In this case, $N$ is an infinitely divisible point process with
canonical measure $\lambda$.
\end{theorem}

In view of (\ref{Kallenberg-condition}), we see that the second
term appearing in the limit of (\ref{cond-L}) is the limit of
$k_n(1-L_{r_n,n}(f))$. Since $n \sim r_n k_n$, the intuition
behind condition (\ref{cond-L}) is that we are forcing
$r_n(L_{m-1,n}(f)-L_{m,n}(f))$, to behave asymptotically as
$1-L_{r_n,n}(f)=\sum_{m=1}^{r_n}(L_{m-1,n}(f)-L_{m,n}(f))$. In
other words, the incremental differences $L_{m-1,n}(f)-L_{m,n}(f)$
with $1 \leq m \leq r_n$, are forced to have the same asymptotic
behavior as their average.

\vspace{3mm}

{\bf Proof:} The proof of the theorem will follow from
(\ref{Kallenberg-condition}), once we show the following relation:
$$\lim_{m\rightarrow m_0}\limsup_{n\rightarrow \infty}
|k_n(1-L_{r_n,n}(f))-n(L_{m-1,n}(f)-L_{m,n}(f))|=0, \quad \forall
f\in C_K^+(E),$$ which can be expressed equivalently as follows,
letting $h(x)=1-e^{-x}$:
\begin{equation}
\label{result-to-prove} \lim_{m\rightarrow
m_0}\limsup_{n\rightarrow \infty} k_n|E[h(N_{r_n,n}(f))]-
r_nE[h(N_{m,n}(f))-h(N_{m-1,n}(f))]|=0, \quad \forall f\in
C_K^+(E).
\end{equation}

In the remaining part of the proof we show that
(\ref{result-to-prove}) holds. Note that only conditions (AD-2)
and (AN), and the stationarity of the array, will be needed for
this. We have
\begin{eqnarray*}
E[h(N_{r_n,n}(f))]&=&E[h(N_{m-1,n}(f))]+\sum_{k=0}^{r_n-m}E[h(N_{k+m,n}(f))-
h(N_{k+m-1,n}(f))] \\
r_nE[h(N_{m,n}(f))-h(N_{m-1,n}(f))]&=&(m-1)E[h(N_{m,n}(f))-h(N_{m-1,n}(f))]+\\
& & \sum_{k=0}^{r_n-m}E[h(N_{k+m,n}(f)-N_{k,n}(f))-
h(N_{k+m-1,n}(f)-N_{k,n}(f))],
\end{eqnarray*}

\noindent where the second equality is due to the strict
stationarity of the sequence $(X_{j,n})_{1 \leq j \leq n}$.
Taking the difference between the previous two equalities, we get:
$$E[h(N_{r_n,n}(f))]-r_nE[h(N_{m,n}(f))-h(N_{m-1,n}(f))]=
mE[h(N_{m-1,n}(f))]-(m-1)E[h(N_{m,n}(f))]
$$
\begin{equation}
\label{interm-step1}
+\sum_{k=0}^{r_n-m}E\{[h(N_{k+m,n}(f))-h(N_{k+m,n}(f)-N_{k,n}(f))]-
[h(N_{k+m-1,n}(f))-h(N_{k+m-1,n}(f)-N_{k,n}(f))] \}.
\end{equation}

\noindent We now apply Taylor's expansion formula:
$h(a)-h(a-b)=b\int_0^1 h'(a-xb)dx$. We get
\begin{eqnarray*}
h(N_{k+m,n}(f))-h(N_{k+m,n}(f)-N_{k,n}(f))&=&N_{k,n}(f)\int_0^1
h'(N_{k+m,n}(f)-x N_{k,n}(f))dx \\
h(N_{k+m-1,n}(f))-h(N_{k+m-1,n}(f)-N_{k,n}(f))&=&N_{k,n}(f)\int_0^1
h'(N_{k+m-1,n}(f)-x N_{k,n}(f))dx.
\end{eqnarray*}

\noindent Taking the difference of the previous two equalities and
applying Taylor's formula again, we obtain:
\begin{eqnarray}
\nonumber
\lefteqn{E|[h(N_{k+m,n}(f))-h(N_{k+m,n}(f)-N_{k,n}(f))]-[h(N_{k+m-1,n}(f))-
h(N_{k+m-1,n}(f)-N_{k,n}(f))] | } \\
\nonumber & & = E \left|N_{k,n}(f)\int_0^1 [h'(N_{k+m,n}(f)-x
N_{k,n}(f))
-h'(N_{k+m-1,n}(f)-x N_{k,n}(f))]dx \right| \\
\nonumber & & = E \left|N_{k,n}(f)(N_{k+m,n}(f)-N_{k+m-1,n}(f))
\int_0^1 h''(\theta_{k,m,n}(x))
 dx \right|=E \left|N_{k,n}(f)f(X_{k+m,n}) \int_0^1
h''(\theta_{k,m,n}(x)) dx \right| \\
\label{interm-step2} & & \leq E [N_{k,n}(f)f(X_{k+m,n}) ],
\end{eqnarray}

\noindent where $\theta_{k,m,n}(x) \geq 0$ is a (random) value
between $N_{k+m-1,n}(f)-x N_{k,n}(f)$ and $N_{k+m,n}(f)-x
N_{k,n}(f)$, and we used the fact that $|h''(\theta)|=e^{-\theta}
\leq 1$ if $\theta \geq 0$. Coming back to (\ref{interm-step1}),
and using (\ref{interm-step2}), we get:
\begin{eqnarray}
\nonumber
k_n|E[h(N_{r_n,n}(f))]-r_nE[h(N_{m,n}(f))-h(N_{m-1,n}(f))]|
&\leq & mk_n E[h(N_{m-1,n}(f))]+(m-1)k_nE[h(N_{m,n}(f))] \\
\label{interm-step3} & & + k_n\sum_{k=0}^{r_n-m} E
[N_{k,n}(f)f(X_{k+m,n}) ].
\end{eqnarray}

We claim that condition (AN) implies:
\begin{equation}
\label{claim1} \lim_{n \to \infty} k_n E[h(N_{m,n}(f))]=0, \quad
\forall m \geq 1.
\end{equation}

\noindent To see this, we use the fact that $h(x) \leq x$ if $x
\geq 0$. If $K$ is the (compact) support of $f$, then
$$k_n E[h(N_{m,n}(f))]  \leq k_n
E\left[\sum_{j=1}^{m}f(X_{j,n})\right]=m k_nE[f(X_{1,n})] \leq m
k_n \|f \|_{\infty}P(X_{1,n} \in K) \leq C \frac{m}{r_n} \|f
\|_{\infty} \to 0.$$

On the other hand, by stationarity,
\begin{eqnarray*}
k_n\sum_{k=0}^{r_n-m} E [N_{k,n}(f)f(X_{k+m,n}) ]&=&
k_n\sum_{k=0}^{r_n-m} \sum_{i=1}^{k}E [f(X_{i,n})f(X_{k+m,n}) \\
&=&k_n \sum_{j=m+1}^{r_n} (r_n-j+1)E [f(X_{1,n})f(X_{j,n})] \leq
n \sum_{j=m+1}^{r_n} E [f(X_{1,n})f(X_{j,n})].
\end{eqnarray*}
Hence, (AD-2) implies that:
\begin{equation}
\label{claim2} \lim_{m \to m_0} \limsup_{n \to \infty}
k_n\sum_{k=0}^{r_n-m} E [N_{k,n}(f)f(X_{k+m,n}) ]=0.
\end{equation}

\noindent Relation (\ref{result-to-prove}) follows from
(\ref{interm-step3}), (\ref{claim1}) and (\ref{claim2}). $\Box$

\vspace{3mm}

The next result shows that if conditions (AD-2) and (\ref{cond-L})
hold with $m_0=1$, then $N$ is a Poisson process.

\begin{proposition}
For each $n \geq 1$, let $(X_{j,n})_{1 \leq j \leq n}$ be a
strictly stationary sequence of $E$-valued random variables.
Suppose that the triangular array $(X_{j,n})_{1 \leq j \leq n, n
\geq 1}$ satisfies condition (AN), as well as conditions (AD-1)
and (AD-2) (with the same sequence $(r_n)_n$). Assume that
$m_0:=\inf {\cS}_{(r_n)_n}=1$, i.e. $$\limsup_{n \to \infty}n
\sum_{j=2}^{r_n} E[f(X_{1,n})f(X_{j,n})]=0, \quad \forall f \in
C_{K}^{+}(E).$$

If there exists a measure $\lambda$ on $M_p(E)\setminus\{o\}$
which satisfies (\ref{cond-lambda-1}), such that
\begin{equation}
\label{cond-L-m0-egal-1} \lim_{n \to \infty}
n(1-E(e^{-f(X_{1,n})})=\int_{M_p(E) \verb2\2 \{o\}}
(1-e^{-\mu(f)})\lambda(d\mu), \quad \forall f\in C_K^+(E),
\end{equation}
then $(N_n)_{n \geq 1}$ converges in distribution to a Poisson
process with intensity $\nu(B):=\lambda(\{\mu \in M_p(E);
\mu(B)=1\})$.
\end{proposition}


{\bf Proof:} By Theorem \ref{main-theorem}, $N_n
\stackrel{d}{\rightarrow} N$, where $N$ is an infinitely divisible
process $N$ with canonical measure $\lambda$.

For each $n \geq 1$, let $(X_{j,n}^{*})_{1 \leq j \leq n}$ be an
i.i.d. sequence with the same distribution as $X_{1,n}$. Let
$N_{n}^{*}=\sum_{j=1}^{n}N_{j,n}^{*}$, where
$N_{j,n}^{*}=\delta_{X_{j,n}^{*}}$. Then $(N_{j,n}^*)_{1 \leq j
\leq n, n \geq 1}$ is a null-array, since
$P(N_{1,n}^{*}(B)>0)=P(X_{1,n} \in B) \to 0$ for all $B \in \cB$.
Note that $(N_{j,n}^{*})_{1 \leq j \leq n}$ are i.i.d. point
processes. By (\ref{cond-L-m0-egal-1}), we have:
$$\lim_{n \to \infty}
\sum_{j=1}^{n}(1-E(e^{-N_{j,n}^{*}(f)})=\int_{M_p(E) \verb2\2
\{o\}} (1-e^{-\mu(f)})\lambda(d\mu), \quad \forall f\in
C_K^+(E).$$

\noindent Therefore, by Theorem 6.1, \cite{kallenberg83}, it
follows that $N_n^* \stackrel{d}{\rightarrow} N$, and $\{n P \circ
[N_{1,n}^*(B_1), \ldots, N_{1,n}^{*}(B_k)]^{-1}\}_n$ converges
weakly to $\lambda \circ \pi_{B_1, \ldots, B_k}^{-1}$, $\forall
B_1, \ldots, B_k \in \cB$. In particular, $nP(X_{1,n} \in
B)=nP(N_{1,n}^{*}(B)=1) \longrightarrow (\lambda \circ
\pi_{B}^{-1})(\{1\})=\nu(B)$, $\forall B \in \cB$, and hence the
sequence $\{n P \circ X_{1,n}^{-1}\}_{n \geq 1}$ converges vaguely
to $\nu$. Since $\lambda$ satisfies (\ref{cond-lambda-1}), the
measure $\nu$ is Radon. By Proposition 3.21, \cite{resnick87}, it
follows that $N_n^* \stackrel{d}{\rightarrow} N^*$, where $N^*$ is
a Poisson process of intensity $\nu$. We conclude that $N
\stackrel{d}{=} N^*$. $\Box$

\section{Partial Sum Convergence}

In this section we suppose that $E=(0,\infty)$. Let
$N_{n}=\sum_{j=1}^{n}\delta_{X_{j,n}}$ and
$S_{n}=\sum_{j=1}^{n}X_{j,n}$.

In Section \ref{weak-convergence}, we have seen various asymptotic
dependence conditions which guarantee the convergence in
distribution of the sequence $(N_n)_{n \geq 1}$ to an infinitely
divisible point process $N$. In the present section, we show that
if the limit process $N$ is ``nice'' (in a sense that will be
specified below), this convergence, together with an asymptotic
negligibility condition in the mean, implies the convergence in
distribution of the partial sum sequence $(S_n)_n$ to an
infinitely divisible random variable. In the literature, this has
been a well-known recipe for obtaining the convergence in
distribution of $(S_n)_n$ to the stable law (see e.g.
\cite{davis83}, \cite{denker-jakubowski89},
\cite{dabrowski-jakubowski94}, \cite{davis-hsing95}). Our
contribution consists in allowing the class of limiting
distributions to include more general infinitely divisible laws.

Let $N$ be an infinitely divisible point process on $(0,\infty)$,
with canonical measure $\lambda$. By Lemma 6.5,
\cite{kallenberg83}, the distribution of $N$ coincides with that
of $\int_{M_p((0,\infty))}\mu \xi(d\mu)$, where $\xi$ is a Poisson
process on $M_p((0,\infty))$ with intensity $\lambda$. Let us
denote by $N_i=\sum_{j \geq 1}\delta_{T_{ij}}, i \geq 1$ the
points of $\xi$, i.e. $\xi=\sum_{i \geq 1}\delta_{N_i}$,
Then the distribution of $N$ coincides with that of $\sum_{i \geq
1}N_i=\sum_{i,j \geq 1}\delta_{T_{ij}}$. (This is called the
``cluster representation'' of $N$.)

The following assumption explains what we meant earlier by a
``nice'' point process $N$.

\begin{assumption}
\label{assumption-lambda} The canonical measure $\lambda$ has the
support contained in the set $M_p^*((0,\infty))$, consisting of
all measures $\mu \in M_{p}((0,\infty))$ whose points are
summable, i.e. all measures $\mu=\sum_{j \geq 1}\delta_{t_j}$ with
$\sum_{j \geq 1}t_j<\infty$.
\end{assumption}

We define the map $T: M_p^*((0,\infty)) \to (0,\infty)$ by
$T(\mu)=\sum_{j \geq 1}t_j$ if $\mu=\sum_{j \geq 1}\delta_{t_j}$.

Assumption \ref{assumption-lambda} is equivalent to saying that
$N_i \in M_p^*((0,\infty))$ a.s. In turn, this is equivalent to
saying that the random variables $U_i:=\sum_{j \geq 1}T_{ij}, i
\geq 1$ are finite a.s. Moreover, we have the following result.

\begin{lemma}
\label{N-is-Poisson} Let $N$ be a point process on $(0,\infty)$
with canonical measure $\lambda$, and the cluster representation:
$$N \stackrel{d}{=} \int_{M_p(E)}\mu \xi(d\mu) = \sum_{i \geq
1}N_i = \sum_{i,j \geq 1}\delta_{T_{ij}}.$$ (Here $\xi=\sum_{i
\geq 1}\delta_{N_i}$ is a Poisson process on $M_p((0,\infty))$
with intensity $\lambda$, and $N_i=\sum_{j \geq 1}\delta_{T_{ij}},
i \geq 1$ are the points of $\xi$.)

Suppose that $\lambda$ satisfies Assumption
\ref{assumption-lambda}, and set $U_{i}:=\sum_{j \geq 1}T_{ij}, i
\geq 1$. Then $N^*:=\sum_{i \geq 1}\delta_{U_i}$ is a Poisson
process with intensity $\rho:=\lambda \circ T^{-1}$, i.e.
$$\rho(A)=\lambda(\{\mu=\sum_{j \geq 1}\delta_{t_j} \in
M_p^*((0,\infty)); \sum_{j \geq 1}t_j \in A \}), \quad \forall A
\in {\cal B}((0,\infty)).$$
\end{lemma}

{\bf Proof:} The lemma will be proved, once we show that for any
measurable $f:(0,\infty) \to (0,\infty)$, we have
$$E\left(e^{-\sum_{i \geq 1}f(U_i)}\right)=\exp\left\{-\int_{0}^{\infty} (1-e^{-f(x)})
\rho(dx)\right\}.$$

\noindent Since $\xi=\sum_{i\geq 1}\delta_{N_i}$ is a Poisson
process with intensity $\lambda$, for any $\psi:M_p((0,\infty))
\to (0,\infty)$ measurable, $$L_{\xi}(\psi)=E\left(e^{-\sum_{i
\geq 1}\psi(N_i)}\right)=\exp\left\{-\int_{M_p^*((0,\infty))}
(1-e^{-\psi(\mu)}) \lambda(d\mu)\right\}.$$
 Let
$\psi_f:M_p^*((0,\infty)) \rightarrow (0,\infty)$ be given by
$\psi_f(\mu)=f(T(\mu))$. Then $\psi_f(N_i)=f(T(N_{i}))=f(\sum_{j
\geq 1}T_{ij})=f(U_i)$ and
$$E\left(e^{-\sum_{i \geq 1}f(U_i)}\right)=\exp\left\{-\int_{M_p^*((0,\infty))}
(1-e^{-\psi_f(\mu)})\lambda(d\mu)\right\}=\exp\left\{-\int_{0}^{\infty}(1-e^{-y})
(\lambda \circ \psi_f^{-1})(dy)\right\}.$$

\noindent  Note that $\rho=\lambda \circ
 T^{-1}$. By the definitions of $\psi_f$ and $\rho$, we have
$\lambda \circ \psi_f^{-1}=\lambda \circ T^{-1} \circ f^{-1}=\rho
\circ f^{-1}$. Hence $E\left(e^{-\sum_{i \geq
1}f(U_i)}\right)=\exp\left\{-\int_{0}^{\infty}(1-e^{-y}) (\rho
\circ
f^{-1})(dy)\right\}=\exp\left\{-\int_{0}^{\infty}(1-e^{-f(x)})
\rho(dx)\right\}$.
 $\Box$

\vspace{3mm}

The next lemma is of general interest and shows that the random
variable $X$ defined as the sum of the points of a Poisson process
on $(0,\infty)$ has an infinitely divisible distribution. To
ensure that $X$ is finite a.s., some restrictions apply to the
intensity $\rho$ of the Poisson process. Recall that a measure
$\rho$ on $(0,\infty)$ is called a {\em L\'evy measure} if
$\int_{(0,1]}x^2 \rho(dx)<\infty$ and $\rho((1,\infty))<\infty$,
or equivalently $\int_0^{\infty}x^2/(1+x^2) \rho(dx)<\infty$.

\begin{lemma}
\label{repr-inf-div-law} Let $N^*=\sum_{i \geq 1}\delta_{U_i}$ be
a Poisson process on $(0,\infty)$, whose intensity $\rho$ is a
L\'evy measure and
\begin{equation}
\label{cond-rho-1} \int_{(0,1]}x\rho(dx)<\infty.
\end{equation}
 Then the random variable
$X:=\sum_{i \geq 1}U_i$ is finite a.s. and has an infinitely
divisible distribution. Moreover,
\begin{equation}
\label{ID-char-function}
E(e^{iuX})=\exp\left\{\int_{0}^{\infty}(e^{iux}-1)\rho(dx)\right\},
\quad \forall u \in \BBr.
\end{equation}
\end{lemma}

{\bf Proof:} Without loss of generality, we can assume that
$U_i=H_{\rho}^{-1}(\Gamma_i)$, where $\Gamma_i=\sum_{j=1}^{i}E_j$,
$(E_j)_{j \geq 1}$ are i.i.d. Exponential$(1)$ random variables,
$H_{\rho}(x)=\rho(x,\infty)$, and
$H_{\rho}^{-1}(y)=\inf\{x>0;H_{\rho}(x) \leq y\}$.

Note that $H_{\rho}$ is a non-increasing function and
$H_{\rho}^{-1}(y) \leq x$ if and only if $y \geq H_{\rho}(x)$.
Then $U_i \leq U_{i-1}, \forall i$ and
\begin{eqnarray}
\nonumber &&P(U_i \leq x_i|U_{1}=x_1, \ldots,
U_{i-1}=x_{i-1})=P(\Gamma_i \geq H_{\rho}(x_i)|
\Gamma_1=H_{\rho}(x_1), \ldots,
\Gamma_{i-1}=H_{\rho}(x_{i-1}))  \\
\nonumber &&=  P(E_i \geq H_{\rho}(x_{i})-H_{\rho}(x_{i-1})|
\Gamma_1=H_{\rho}(x_1),
\ldots, \Gamma_{i-1}=H_{\rho}(x_{i-1})) \\
\label{allow-FK} && = P(E_i \geq
H_{\rho}(x_{i})-H_{\rho}(x_{i-1}))
=e^{-(H_{\rho}(x_{i})-H_{\rho}(x_{i-1}))} \ \quad \mbox{for all} \
x_i \leq x_{i-1} \leq \ldots \leq x_1
\end{eqnarray}

\noindent Relation (\ref{allow-FK}) allows us to invoke a powerful
(and highly non-trivial) construction, due to Ferguson and Klass
(see \cite{ferguson-klass72}). More precisely, let $(V_i)_{i \geq
1}$ be a sequence of i.i.d. random variables with values in
$[0,1]$ and common distribution $G$, which is independent of
$(U_i)_{i \geq 1}$, and define $Y_{t}=\sum_{i \geq 1}U_i 1_{\{V_i
\leq t\}}, t \in [0,1]$. Then, Ferguson and Klass showed that
$(Y_t)_{t \in [0,1]}$ is a L\'evy process with characteristic
function
$E(e^{iuY_t})=\exp\left\{G(t)\int_{0}^{\infty}(e^{iux}-1)\rho(dx)\right\}$,
$\forall u \in \BBr$. The proof is complete by observing that
$X=Y_1=\sum_{i \geq 1}U_i$. $\Box$




\begin{example}
{\rm $\rho(dx)=\alpha x^{-1}e^{-x}1_{\{x>0\}}dx$ with $\alpha>0$.
In this case, $X$ has a Gamma($\alpha$) distribution.}
\end{example}

\begin{example}
{\rm $\rho(dx)=c_{\alpha} x^{-\alpha-1}1_{\{x>0\}}dx$ with $\alpha
\in (0,1)$. In this case, $X$ has a stable distribution of index
$\alpha$.}
\end{example}

As a by-product of the previous lemma, we obtain a representation
of an infinitely divisible distribution, similar to the
LePage-Woodroofe-Zinn representation of the stable law (Theorem 2,
\cite{LWZ81}). The proof of this corollary is based on a
representation of a Poisson process, which is included in Appendix
A.

\begin{corollary}
Let $\rho$ be a measure on $(0,\infty)$, which is given by the
following ``product-convolution'' type formula:
\begin{equation}
\label{cond-rho} \rho(A)=\int_{0}^{\infty}
\int_{0}^{\infty}1_{A}(wy)F(dw)\nu(dy), \quad \forall A \in {\cal
B}((0,\infty)),
\end{equation}
where $\nu$ is an arbitrary Radon measure $\nu$ on $(0,\infty)$
and $F$ is an arbitrary probability measure on $(0,\infty)$.

If the measure $\rho$ is L\'evy and satisfies (\ref{cond-rho-1}),
then any infinitely divisible random variable $X$ with
characteristic function (\ref{ID-char-function}) admits the
representation $X\stackrel{d}{=}\sum_{i \geq 1}P_iW_i$, where
$(P_i)_{i \geq 1}$ are the points of a Poisson process of
intensity $\nu$, and $(W_i)_{i \geq 1}$ is an independent i.i.d.
sequence with distribution $F$.
\end{corollary}

{\bf Proof:} Let $N^*$ be a Poisson process on $(0,\infty)$, of
intensity $\rho$. Using definition (\ref{cond-rho}) of $\rho$, and
by invoking Proposition \ref{productPiWi} (Appendix A) , it
follows that $N^*$ admits the representation $N^*
\stackrel{d}{=}\sum_{i \geq 1}\delta_{P_iW_i}$, where $(P_i)_{i
\geq 1}$  and $(W_i)_{i \geq 1}$ are as in the statement of the
corollary. By Lemma \ref{repr-inf-div-law}, it follows that the
random variable $X:=\sum_{i \geq 1}P_iW_i$ has an infinitely
divisible distribution with characteristic function
(\ref{ID-char-function}). $\Box$

\begin{remark}
\label{stable-case-calculations} {\rm If we let $\nu(dx)=\alpha
x^{-\alpha-1}1_{\{x>0\}}dx$ and $F$ be an arbitrary probability
measure $F$ on $(0,\infty)$, then the measure $\rho$ given by
(\ref{cond-rho}) satisfies:
$$\rho(x,\infty)=
\int_{0}^{\infty}
  \int_{0}^{\infty} 1_{(x,\infty)} (wy) F(dw) \nu(dy)
=\int_{0}^{\infty} \int_{0}^{\infty} \nu\left(\frac{x}{w},\infty
\right)F(dw)= x^{-\alpha}\int_{0}^{\infty}
w^{\alpha}F(dw)=\gamma_{\alpha} x^{-\alpha},$$ where
$\gamma_{\alpha}=\int_0^{\infty}w^{\alpha}F(dw)$. Hence
$\rho(dx)=\alpha \gamma_{\alpha}x^{-\alpha-1}1_{\{x>0\}}dx$.}
\end{remark}

To obtain the convergence of the partial sum sequence, we
introduce a new asymptotic negligibility condition.

\begin{definition}
$(X_{j,n})_{1 \leq j \leq n, n \geq 1}$ satisfies {\bf condition
(AN')} if $\lim_{\varepsilon \rightarrow 0} \limsup_{n \rightarrow
\infty} nE[X_{1,n}1_{\{X_{1,n} \leq \varepsilon\}}]=0$.
\end{definition}

The next theorem is a generalization of Theorem 3.1,
\cite{davis-hsing95}, to the case of an arbitrary infinitely
divisible law (without Gaussian component, and whose L\'evy
measure $\rho$ satisfies (\ref{cond-rho-1})), as the limiting
distribution of $(S_n)_n$.

\begin{theorem}
\label{partial-sum-theorem} For each $n \geq 1$, let $(X_{j,n})_{1
\leq j \leq n}$ be a strictly stationary sequence of positive
random variables. Suppose that the array $(X_{j,n})_{1 \leq j \leq
n, n \geq 1}$ satisfies condition (AN'). Let
$N_n=\sum_{j=1}^{n}\delta_{X_{j,n}}$ and
$S_n=\sum_{j=1}^{n}X_{j,n}$.

If

 (i) $N_n \stackrel{d}{\rightarrow} N$, where $N$ is an infinitely
divisible point process, whose canonical measure $\lambda$
satisfies Assumption \ref{assumption-lambda}; and

(ii) $\rho:=\lambda \circ T^{-1}$ is a L\'evy measure and
satisfies (\ref{cond-rho-1}),

\noindent then $(S_n)_n$ converges in distribution to an
infinitely divisible random variable with characteristic function
(\ref{ID-char-function}).
\end{theorem}

{\bf Proof:} For each $\varepsilon>0$ arbitrary, we write
\begin{equation}
\label{Sn-decomposition-eta} S_n=S_{n}(\varepsilon,
\infty)+S_{n}(0,\varepsilon).
\end{equation}
where $S_{n}(\varepsilon,\infty)=\sum_{j=1}^{n}X_{j,n} 1_{\{
X_{j,n}>\varepsilon\}}$ and
$S_{n}(0,\varepsilon)=\sum_{j=1}^{n}X_{j,n}1_{\{X_{j,n}
 \leq \varepsilon\}}$. Let $N \stackrel{d}{=}\sum_{i,j \geq
 1}\delta_{T_{ij}}$ be the cluster representation of $N$ and
 $U_i=\sum_{j \geq 1}T_{ij}$ for all $i \geq 1$.

By Lemma \ref{N-is-Poisson}, $N^*:=\sum_{i \geq 1}\delta_{U_i}$ is
a Poisson process of intensity $\rho:=\lambda \circ T^{-1}$. By
Lemma \ref{repr-inf-div-law}, the random variable $X:=\sum_{i \geq
1}U_i=\sum_{i,j \geq 1}T_{ij}$ is finite a.s. and has an
infinitely divisible distribution. Moreover,
(\ref{ID-char-function}) holds.

Define $T_{\varepsilon}: M_p^*((0,\infty)) \rightarrow (0,\infty)$
by $T_{\varepsilon}(\mu=\sum_{j \geq 1}\delta_{t_j})=\sum_{j \geq
1} t_j 1_{\{t_{j}> \varepsilon\}}$.  Note that $T_{\varepsilon}$
is continuous $P \circ N^{-1}$-a.s. By the continuous mapping
theorem, we get $T_{\varepsilon}(N_n)=S_n(\varepsilon, \infty)
\stackrel{d}{\rightarrow} T_{\varepsilon}(N)=\sum_{i,j \geq
1}T_{ij} 1_{\{T_{ij}> \varepsilon\}}$, as $n \to \infty$. Since
$X=\sum_{i,j \geq 1}T_{ij}$ converges a.s., it follows that
$\sum_{i,j \geq 1}T_{ij} 1_{\{T_{ij} > \varepsilon\}}
\stackrel{a.s}{\to} X=\sum_{i,j \geq 1}T_{ij}$ as $\varepsilon \to
0$. Hence
\begin{equation}
\label{Sn-upper-epsilon} S_n(\varepsilon,\infty) \stackrel{d}{\to}
X \ \quad {\rm as} \ n \to \infty, \varepsilon \to 0.
\end{equation}

\noindent By Markov's inequality and condition (AN'), we see that
for any $\delta>0$, $P(S_{n}(0,\varepsilon)>\delta) \leq
\delta^{-1} E[S_{n}(0,\varepsilon)]=\delta^{-1} n
E[X_{1,n}1_{\{X_{1,n} \leq \varepsilon\}}] \rightarrow 0$, as $n
\to \infty, \varepsilon \to 0$. Hence
\begin{equation}
\label{Sn-lower-epsilon} S_{n}(0,\varepsilon) \stackrel{P}{\to}0
\quad \mbox{as} \ n \to \infty, \varepsilon \to 0.
\end{equation}

\noindent From (\ref{Sn-decomposition-eta}),
(\ref{Sn-upper-epsilon}) and (\ref{Sn-lower-epsilon}), we conclude
that $S_{n} \stackrel{d}{\rightarrow} X$.

$\Box$

\begin{remark}
{\rm Lemma \ref{repr-inf-div-law} can be extended to a Poisson
process whose intensity $\rho$ is an arbitrary L\'evy measure on
$(0,\infty)$. More precisely, using the Theorem of
\cite{ferguson-klass72}, one can prove that if $N^*=\sum_{i \geq
1}\delta_{U_i}$ is a Poisson process on $(0,\infty)$, whose
intensity $\rho$ is a L\'{e}vy measure,
then the random variable $Y:=\sum_{i \geq 1}(U_i-c_i)$ is finite
a.s. and has an infinitely divisible distribution. Moreover,
$$E(e^{iuY})=\exp\left\{\int_{0}^{\infty}\left(e^{iux}-1-
\frac{iux}{1+x^2}\right) \rho(dx)\right\}, \quad \forall u \in
\BBr,$$ where the constants $c_i$ are defined by:
$c_i=\int_{H_{\rho}^{-1}(i)}^{H_{\rho}^{-1}(i-1)}x/(1+x^2)\rho(dx)$.
If $\gamma=\sum_{i \geq 1}c_i=\int_{0}^{\infty}x/(1+x^2)\rho(dx)$
is finite, then one can conclude that the random variable
$X=\sum_{i \geq 1}U_i=Y+\gamma$ (which appears in Theorem
\ref{partial-sum-theorem}) has an infinitely divisible
distribution with characteristic function
$$E(e^{iuX})=\exp\left\{iu\gamma +\int_{0}^{\infty}\left(e^{iux}-1-
\frac{iux}{1+x^2}\right) \rho(dx)\right\}, \quad \forall u \in
\BBr.$$ Unfortunately, requiring that $\gamma$ is finite is
equivalent to saying that $\int_{(0,1]}x \rho(dx)<\infty$, which
is precisely the restriction imposed on $\rho$ in Lemma
\ref{repr-inf-div-law}. In other words, condition
(\ref{cond-rho-1}) cannot be removed from Theorem
\ref{partial-sum-theorem}, using the Ferguson and Klass approach.}
\end{remark}

We finish this section with an example for which the hypothesis of
Theorem \ref{partial-sum-theorem} are verified. This example is
based on the recent work \cite{kulik06}, generalizing the moving
average model MA$(\infty)$ to the case of random coefficients.

\begin{example}
{\rm (Linear processes with random coefficients) Let
$X_{i,n}=X_i/a_n$ for all $1 \leq i \leq n$, where
$$X_{i}=\sum_{j=0}^{\infty}C_{i,j}Z_{i-j} \quad \mbox{for all} \ i \geq 1.$$
The objects $(Z_k)_{k \in \BBz}$, $(a_n)_{n \geq 1}$ and
$(C_{i,j})_{i \geq 1,j \geq 0}$ are defined as follows:
\begin{itemize}

\item $(Z_{k})_{k \in \BBz}$ is a sequence of i.i.d. positive
random variables such that $Z_0 \stackrel{d}{=}Z$, where $Z$ has
heavy tails, i.e. $P(Z>x)=x^{-\alpha}L(x)$ for $\alpha \in (0,2)$
and $L$ a slowly varying function.

\item $(a_n)_{n \geq 1}$ is a non-decreasing sequence of positive
numbers such that $P(Z>a_n) \sim n^{-1}$.

\item $(C_{i,j})_{i \geq 1,j \geq 0}$ is an array of positive
random variables, which are independent of $(Z_k)_{k \in \BBz}$.
We suppose that the rows  $(C_{1,j})_{j \geq 0}, (C_{2,j})_{j \geq
0}, \ldots$ of this array are i.i.d. copies of a sequence
$(C_{j})_{j \geq 0}$ of positive random variables.
Moreover, we suppose that the sequence $(C_j)_{j \geq 0}$
satisfies certain moment conditions, which imply that
$c:=\sum_{j=0}^{\infty}E[C_j^{\alpha}]<\infty$. (We refer the
reader to condition (D) of \cite{kulik06} for the exact moment
conditions. In fact, we may allow for a mixing-type dependence
structure between the rows.)

\end{itemize}

Proposition 2.1, \cite{kulik06} shows that $P(X_1>x) \sim c
P(Z>x)$ as $x \to \infty$. Since $Z$ has heavy tails, it follows
that $X_1$ has heavy tails too. Assume that $\alpha \in (0,1)$.
In this case, one can prove that: (see e.g. (3.6) in
\cite{davis-hsing95})
$$\lim_{\varepsilon \to 0} \limsup_{n \to \infty} \frac{n}{a_n}E[X_{1}1_{\{X_{1,n}
\leq a_n \varepsilon\}}]=0,$$ i.e. the array $(X_{j,n})_{1 \leq j
\leq n, n \geq 1}$ satisfies condition (AN').

Let $N_n=\sum_{i=1}^{n}\delta_{X_{i,n}}$. By Theorem 3.1,
\cite{kulik06}, $N_{n} \stackrel{d}{\rightarrow} N$, where $N$ is
an infinitely divisible point process with the cluster
representation $N \stackrel{d}{=}\sum_{i \geq 1} \sum_{j \geq
0}\delta_{P_i C_{i,j}}$. Here $(P_i)_{i \geq 1}$ are the points of
a Poisson process of intensity $\nu(dx)=\alpha
x^{-\alpha-1}1_{\{x>0\}}dx$, which is independent of the array
$(C_{i,j})_{i \geq 1, j \geq 0}$. Since $\alpha \in (0,1)$, it
follows that $W_i:=\sum_{j \geq 0}C_{i,j}<\infty$ a.s. for all $i
\geq 1$. Hence, the random variables $U_i:=\sum_{j \geq 0}P_i
C_{i,j}=P_iW_i, i \geq 1$ are finite a.s. and Assumption
\ref{assumption-lambda} is verified. This proves that condition
(i) in Theorem \ref{partial-sum-theorem} is satisfied.

By Lemma \ref{N-is-Poisson}, the process $N^*:=\sum_{i \geq
1}\delta_{U_i}=\sum_{i \geq 1}\delta_{P_iW_i}$ is a Poisson
process of intensity $\rho:=\lambda \circ T^{-1}$, where $\lambda$
is the canonical measure of $N$. From Proposition
\ref{productPiWi} (Appendix A), it follows that $\rho$ satisfies
(\ref{cond-rho}). By Remark \ref{stable-case-calculations}, it
follows that $\rho(dx)=\alpha
\gamma_{\alpha}x^{-\alpha-1}1_{\{x>0\}}$, where
$\gamma_{\alpha}=\int_0^{\infty}w^{\alpha}F(dw)$. Clearly, this
measure $\rho$ is L\'evy; it satisfies condition
(\ref{cond-rho-1}) since $\alpha<1$. This proves that condition
(ii) in Theorem \ref{partial-sum-theorem} is satisfied.

By applying Theorem \ref{partial-sum-theorem}, it follows that
$(S_n)_{n \geq 1}$ converges in distribution to an infinitely
divisible law with characteristic function
(\ref{ID-char-function}), which is in fact the stable law of index
$\alpha$.}

\end{example}

\section{Real Valued Observations}

In this section, we assume that $E=\BBr \verb2\2 \{0\}$. By Lemma
2.1, \cite{davis-hsing95}, the support of the canonical measure
$\lambda$ of an infinitely divisible point process on $\BBr
\verb2\2 \{0\}$, is contained in the set $M_0(\BBr \verb2\2
\{0\})$, defined by:
$$M_0(\BBr \verb2\2 \{0\})=\{\mu=\sum_{j
\geq 1} \delta_{t_j} \in M_p(\BBr \verb2\2 \{0\}) \verb2\2 \{o\};
\ \exists \ x_{\mu} \in (0,\infty) \ \mbox{such that} \ |t_j| \leq
x_\mu \forall j \geq 1\}.$$
 Let $\tilde M(\BBr \verb2\2 \{0\})=\{\mu
\in M_0(\BBr \verb2\2 \{0\}); |t_j| \leq 1, \forall j \geq 1 \}$.
The following result gives the necessary and sufficient condition
for a ``product-type'' cluster representation of an infinitely
divisible process, as in Corollary 2.4, \cite{davis-hsing95}.

\begin{proposition}
\label{representation-ID-process} Let $N$ be an infinitely
divisible point process on $\BBr \verb2\2 \{0\}$, with canonical
measure $\lambda$. Then $N\stackrel{d}{=}\sum_{i,j \geq
1}\delta_{P_i Q_{ij}}$,  where $(P_i)_{i \geq 1}$ are the points
of a Poisson process on $(0,\infty)$ of intensity $\nu$ and
$(Q_{1j})_{j \geq 1}, (Q_{2j})_{j \geq 1}, \ldots$ are i.i.d.
sequences, independent of $(P_i)_{i \geq 1}$, if and only if there
exists a probability measure ${\cal O}$ on $\tilde M(\BBr \verb2\2
\{0\})$ such that, for every measurable non-negative function $f$
on $\BBr \verb2\2 \{0\}$, we have
$$\int_{M_0(\BBr \verb2\2 \{0\})} (1-e^{-\mu(f)})\lambda(d\mu)
=\int_0^{\infty} \int_{\tilde M(\BBr \verb2\2
\{0\})}(1-e^{-\mu(f(y \cdot))}){\cal O}(d \mu) \nu(dy).$$ In this
case, ${\cal O}$ is the distribution of $\sum_{j \geq
1}\delta_{Q_{1j}}$.
\end{proposition}

{\bf Proof:} Let $N'=\sum_{i,j \geq 1}\delta_{P_{i}Q_{ij}}$.
Clearly,
 $L_{N}(f)=\exp\left\{-\int_{M_0(\BBr \verb2\2 \{0\})}
(1-e^{-\mu(f)})\lambda(d\mu) \right\}$. Following the same lines
as for the proof of (\ref{Laplace-N**}) (Appendix A), one can show
that
$$L_{N'}(f)=\exp\left\{-\int_0^{\infty}\int_{\tilde M(\BBr \verb2\2
\{0\})}(1-e^{-\mu(f(y \cdot))}){\cal O}(d \mu) \nu(dy)\right\}.$$
The result follows since $N \stackrel{d}{=} N'$ if and only if
$L_N(f)=L_{N'}(f)$ for every measurable function $f$. $\Box$

\vspace{3mm}

In the light of Proposition \ref{representation-ID-process}, the
following result becomes an immediate consequence of Theorem
\ref{main-theorem}.

\begin{corollary}
\label{corollary-main-theorem} For each $n \geq 1$, let
$(X_{j,n})_{1 \leq j \leq n}$ be a strictly stationary sequence of
random variables with values in  $\BBr \verb2\2 \{0\}$. Suppose
that the array $(X_{j,n})_{1 \leq j \leq n, n \geq 1}$ satisfies
condition (AN), as well as conditions (AD-1) and (AD-2) (with the
same sequence $(r_n)_n$). Let $m_0:=\inf {\cS}_{(r_n)_n}$.

If there exists a Radon measure $\nu$ on $(0,\infty)$ and a
probability measure ${\cal O}$ on $\tilde M(\BBr \verb2\2 \{0\})$,
such that
\begin{equation}
\label{E1T1} \lim_{m\rightarrow m_0}\limsup_{n\rightarrow \infty}
\left|n\left(L_{m-1,n}(f)-L_{m,n}(f)\right)-\int_0^{\infty}\int_{\tilde
M(\BBr \verb2\2 \{0\})}(1- e^{-\mu f(y\cdot)}) {\cal O}(d\mu)
\nu(dy)\right|=0, \quad \forall f \in C_{K}^{+}(\BBr \verb2\2
\{0\}),
\end{equation}
then $N_n \stackrel{d}{\rightarrow} N$, where $N=\sum_{i,j \geq
1}\delta_{P_iQ_{ij}}$, $(P_i)_{i \geq 1}$ are the points of a
Poisson process on $(0,\infty)$ of intensity $\nu$, and
$(Q_{1j})_{j \geq 1}, (Q_{2j})_{j \geq 1}, \ldots$ are i.i.d.
sequences with distribution ${\cal O}$, independent of $(P_i)_{i
\geq 1}$.
\end{corollary}

\begin{remark}
{\rm In particular, one may restate Corollary
\ref{corollary-main-theorem}, in the case $X_{j,n}=X_{j}/a_n$,
where $(X_j)_{j \geq 1}$ is a strictly stationary sequence of
random variables with values in $\BBr \verb2\2 \{0\}$ such that
$X_1$ has heavy tails, and $(a_n)_{n \geq 1}$ satisfies
$P(X_1>a_n) \sim n^{-1}$. The result obtained in this manner can
be viewed as a complement to Theorem 2.3, \cite{davis-hsing95}.}
\end{remark}

Recall that a bounded Borel set in $\BBr \verb2\2 \{0\}$ is
bounded away from $0$. Therefore, condition (AN) holds if
$\limsup_{n\rightarrow \infty}nP(|X_{1,n}|\geq \eta )<\infty$,
$\forall \eta>0$. Note also that condition (AD-2') holds if there
exists a sequence $(r_n)_n \subset \BBz_{+}$ with $r_n \to \infty$
such that
\begin{equation}
\label{AD-2''} \lim_{m\rightarrow m_0}\limsup_{n\rightarrow
\infty}n\sum_{j=m+1}^{r_n} P(|X_{1,n}|\geq \eta, |X_{j,n}|\geq
\eta)=0, \quad \forall \eta>0,
\end{equation}
where $m_0:=\inf {\cS}_{(r_n)_n}$.
Condition (\ref{AD-2''}) 
can be viewed as an asymptotic``anti-clustering'' condition for
the process $S_{n}(t)=\sum_{j=1}^{[nt]}X_{j,n}, t \in [0,1]$. To
see this, note that this cadlag process jumps at times $t_j=j/n$
with $1 \leq j \leq n$, the respective jump heights being $\Delta
S_{n}(t_j)=X_{j,n}$. By (\ref{AD-2-stationarity}), we have
$$P(\exists i<k\leq r_n \ \mbox{with} \ k-i \geq m \
 \mbox{such that} \ |\Delta S_n(t_i) |\geq \eta,
|\Delta S_{n}(t_k)|\geq \eta ) \leq r_n\sum_{j=m+1}^{r_n}
P(|X_{1,n}|\geq \eta, |X_{j,n}|\geq \eta ).$$

Therefore, we can explain intuitively condition (\ref{AD-2''}) by
saying that the chance that the process $(S_n(t))_{t\in[0,1]}$ has
at least two jumps that exceed $\eta$ in the time interval $[0,
r_n/n]$ (and are located at a minimum distance of $m/n$ of each
other) is asymptotically zero. (See also p. 213, \cite{ekm97}.)


\section{Examples of arrays satisfying (AD-1)}
\label{section-about-AD}

In this section, we assume that $E=\BBr \verb2\2 \{0\}$ and we
examine condition (AD-1) in the case of arrays which possess a
known dependence structure on each row.

\subsection{$m$-dependent or strongly mixing sequences}

Recall that the $m$-th order mixing coefficient of a sequence
$(X_{j})_{j \geq 1}$ of random variables is defined by:

$$\alpha(m)=\sup\{|P(A \cap B)-P(A)P(B)|; A \in \sigma(X_1, \ldots, X_k), B \in \sigma(X_{k+m},
X_{k+m+1}, \ldots), k \geq 1 \}.$$ The random variables $(X_j)_{j
\geq 1}$ are called {\em strongly mixing} if $\lim_{m \to
\infty}\alpha(m) = 0$.

If $X$ is a $\sigma(X_1, \ldots, X_k)$-measurable bounded random
variable and $Y$ is a $\sigma(X_{k+m},
X_{k+m+1},\ldots)$-measurable bounded random variable, then: (see
e.g. \cite{ibragimov62})
\begin{equation}
\label{Bradley-Bryc-inequality} |E(XY)-E(X)E(Y)| \leq 4 \alpha(m)
\|X\|_{\infty} \|Y \|_{\infty}.
\end{equation}

\begin{lemma}
\label{AD1-m-dependent} For each $n \geq 1$, let $(X_{j,n})_{1
\leq j \leq n}$ be a strictly stationary sequence of random
variables and $\alpha_n(m)$ be its $m$-th order mixing
coefficient, for $m<n$. Suppose that either

(i) $\alpha_n(m')=0$ for all $m' \geq m$, $n \geq 1$; or

(ii) $\alpha_n(m)=\alpha(m) \quad  \forall n> m, \forall m \geq 1$
and $\lim_{m \to \infty}\alpha(m) =0$.

\noindent If the triangular array $(X_{j,n})_{1 \leq j \leq n, n
\geq 1}$ satisfies condition (AN), then it also satisfies
condition (AD-1).
\end{lemma}

\begin{remark}
{\rm a) Condition (i) requires that the sequence $(X_{j,n})_{1
\leq j \leq n}$ is $m$-dependent, for any $n \geq 1$.

\noindent b) Condition (ii) is satisfied if $X_{j,n}=X_{j}/a_n$
and $(X_j)_{j \geq 1}$ a strictly stationary strongly mixing
sequence. }
\end{remark}

{\bf Proof:}  We want to prove that there exists a sequence
$(r_n)_n \to \infty$ with $k_n:=[n/r_n] \to \infty$ such that
\begin{equation}
\label{AD1-m-dep-step1} E
(e^{-N_{n}(f)})-\{E(e^{-N_{r_n,n}(f)})\}^{k_n} \to 0, \quad
\forall f \in C_{K}^{+}(E).
\end{equation}

\noindent Note that
$e^{-N_{k_nr_n,n}(f)}-e^{-N_{n}(f)}=e^{-N_{k_nr_n,n}(f)}
(1-e^{-\sum_{j=k_nr_n+1}^{n}f(X_{j,n})}) \leq
\sum_{j=k_nr_n+1}^{n}f(X_{j,n})$, using the fact that $1-e^{-x}
\leq x$ for any $x \geq 0$. By stationarity and condition (AN), we
obtain that:
$$ \left|E ( e^{-N_{n}(f)})-E(e^{-N_{k_nr_n,n}(f)})\right| \leq
(n-r_nk_n)E[f(X_{1,n})]\leq r_n \|f\|_{\infty} P(X_{1,n}\in K)
\leq \frac{1}{k_n} \|f\|_{\infty} C \to 0,
$$
where $K$ is the compact support of $f$. Therefore, in order to
prove (\ref{AD1-m-dep-step1}), it is enough to show that
\begin{equation}
\label{AD1-mdep} E
(e^{-N_{k_nr_n,n}(f)})-\{E(e^{-N_{r_n,n}(f)})\}^{k_n} \to 0.
\end{equation}

To prove (\ref{AD1-mdep}), we will implement Jakubowski's ``block
separation'' technique (see the proof of Proposition 5.2,
\cite{jakubowski93}, for a variant of this technique). Let
$(m_n)_n$ be a sequence of positive integers such that
\begin{equation}
\label{properties-of-m-n}
 m_n \to \infty, \quad m_n/r_n \to 0, \quad
\mbox{and} \quad k_n \alpha(m_n) \to 0.
\end{equation}
(The construction of sequences $(r_n)_n$ and $(m_n)_n$ which
satisfy (\ref{properties-of-m-n}) is given in Appendix B.)

For each $n \geq 1$, we consider $k_n$ blocks of consecutive
integers of length $r_n-m_n$, separated by ``small'' blocks of
length $m_n$:
\begin{center}
\begin{picture}(400,20)

\put(0,10){\line(1,0){340}}

\put(0,0){$0$}    \put(0,10){\line(0,1){3}}

\put(98,7){$\times$} \put(80,10){\line(0,1){3}} \put(97,0){$r_n$}
\put(60,0){$r_n-m_n$} 

\put(198,7){$\times$} \put(180,10){\line(0,1){3}}
\put(195,0){$2r_n$}
\put(153,0){$2r_n-m_n$} 

\put(250,0){$\ldots$}

\put(317,10){\line(0,1){3}} \put(335,0){$k_nr_n$}
\put(337,7){$\times$} \put(285,0){$k_n r_n-m_n$}

\end{picture}
\end{center}

More precisely, for each $1 \leq i \leq k_n$, let $H_{i,n}$ be the
(big) block of consecutive integers between $(i-1)r_n+1$ and
$ir_n-m_n$ and $I_{i,n}$ be the (small) block of size $m_n$,
consisting of the integers between $ir_n-m_n+1$ and $ir_n$. Let
$$U_{i,n}=\sum_{j \in H_{i,n}}f(X_{j,n})=N_{ir_n-m_n,n}(f)-N_{(i-1)r_n,n}(f).$$

\noindent By the stationarity of the array, $(U_{i,n})_{1 \leq i
\leq k_n}$ are identically distributed. Clearly,
$U_{1,n}=N_{r_n-m_n,n}(f)$.

On the other hand, since the separation blocks have size $m_n$,
which is ``relatively small'' compared to $r_n$,
\begin{eqnarray}
\label{AD-first-conv} \lim_{n \to \infty}
|E(e^{-N_{k_nr_n,n}(f)})-E(e^{-\sum_{i=1}^{k_n}U_{i,n}})| &= &0 \\
\label{AD-second-conv} \lim_{n \to \infty}
|\{E(e^{-N_{r_n,n}(f)})\}^{k_n}-\{E(e^{-U_{1,n}})\}^{k_n}|& =& 0.
\end{eqnarray}

\noindent (To prove (\ref{AD-first-conv}),  note that
$e^{-\sum_{i=1}^{k_n}U_{i,n}}-e^{-N_{k_nr_n,n}(f)}=
e^{-\sum_{i=1}^{k_n}U_{i,n}}(1-e^{-\sum_{i=1}^{k_n}\sum_{j \in
I_{i,n}}f(X_{j,n})}) \leq \linebreak \sum_{i=1}^{k_n}\sum_{j \in
I_{i,n}}f(X_{j,n})$, using the fact that $1-e^{-x} \leq x$ for any
$x \geq 0$. Hence
$|E(e^{-N_{k_nr_n,n}(f)})-E(e^{-\sum_{i=1}^{k_n}U_{i,n}^{(m)}})|
\leq m_nk_n E[f(X_{1,n})] \leq m_n k_n \|f\|_{\infty} P(X_{1,n}
\in K) \leq   (m_n/r_n) \|f\|_{\infty} C \to 0$, where we used
condition (AN) and (\ref{properties-of-m-n}). Relation
(\ref{AD-second-conv}) follows by a similar argument, using the
fact that $|x^k-y^k| \leq k|x-y|$ for any $x,y \geq 0$ and $k \in
\BBz_{+}$.)

Therefore, in order to prove (\ref{AD1-mdep}), it suffices to show
that:
\begin{equation}
\label{AD-third-conv} \lim_{n \to \infty}
|E(e^{-\sum_{i=1}^{k_n}U_{i,n}})-\{E(e^{-U_{1,n}})\}^{k_n}| =0.
\end{equation}

In case (i), this follows immediately since the random variables
$(U_{i,n})_{1 \leq i \leq k_n}$ are independent, for $n$ large.

In case (ii), we claim that, for any $1 \leq k \leq k_n$ we have
\begin{equation}
\label{induction}
|E(e^{-\sum_{i=1}^{k}U_{i,n}})-\{E(e^{-U_{1,n}})\}^{k}| \leq 4
(k-1) \alpha(m_n).
\end{equation}

\noindent (Relation (\ref{induction}) can be proved by induction
on the number $k$ of terms. If $k=2$, then (\ref{induction})
follows from inequality (\ref{Bradley-Bryc-inequality}). If
relation (\ref{induction}) holds for $k-1$, then
$|E(e^{-\sum_{i=1}^{k}U_{i,n}})-\{E(e^{-U_{1,n}})\}^{k}| \leq
|E(e^{-\sum_{i=1}^{k}U_{i,n}})-
E(e^{-\sum_{i=1}^{k-1}U_{i,n}})E(e^{-U_{k,n}})|
+|E(e^{-\sum_{i=1}^{k-1}U_{i,n}})- \{E(e^{-U_{1,n}})\}^{k-1}| $.
For the first term we use (\ref{Bradley-Bryc-inequality}), since
$\sum_{i=1}^{k-1}U_{i,n}$ and $U_{k,n}$ are separated by a block
of length $m_n$. For the second term we use the induction
hypothesis.)

From (\ref{induction}) and (\ref{properties-of-m-n}), we get:
$$|E(e^{-\sum_{i=1}^{k_n}U_{i,n}})-\{E(e^{-U_{1,n}})\}^{k_n}| \leq
4 k_n \alpha(m_n) \to 0.$$
$\Box$

\subsection{Associated sequences}

Recall that the random variables $(X_{j})_{j \geq 1}$ are called
{\em associated} if for any finite disjoint sets $A,B$ in $\{1,2,
\ldots \}$ and for any coordinate-wise non-decreasing functions
$h: \BBr^{\# A} \to \BBr$ and $k : \BBr^{\# B} \to \BBr$
$${\rm Cov}(h(X_{j},j \in A), k(X_{j},j \in B)) \geq 0,$$
where $\# A$ denotes the cardinality of the set $A$. (See e.g.
\cite{barlow-prochan81}, \cite{epw67} for more details about the
association.)

If $(X_{j})_{j \geq 1}$ is a sequence of associated random
variables, then for any finite disjoint sets $A,B$ in $\{1,2,
\ldots \}$ and for any functions $h: \BBr^{\# A} \to \BBr$ and $k
: \BBr^{\# B} \to \BBr$ (not necessarily coordinate-wise
non-decreasing), which are partially differentiable and have
bounded partial derivatives, we have: (see Lemma 3.1.(i),
\cite{birkel88})
\begin{equation}
\label{birkel-inequality} |{\rm Cov}(h(X_{j},j \in A), k(X_{j},j
\in B)) | \leq \sum_{i \in A} \sum_{j \in B} \left\|\frac{\partial
h}{\partial x_i} \right\|_{\infty}
 \left\|\frac{\partial k}{\partial x_j} \right\|_{\infty}
 {\rm Cov}(X_{i},X_{j}).
 \end{equation}

Let $\cC$ be the class of all bounded nondecreasing functions $g$,
for which there exists a compact subset $K$ of $E$ such that
$g(x)=x$ for all $x \in K$. Let ${\cal S}_1$ be the set of all $m
\in \BBz_{+}$ for which $$\limsup_{n \to \infty}n
\sum_{j=m+1}^{n}{\rm Cov}(g(X_{1,n}),g(X_{j,n}))=0 \quad \forall \
g \in \cC.$$ Let $m_1$ be the smallest integer in ${\cal S}_1$. By
convention, we let $m_1=\infty$ if ${\cal S}_1=\emptyset$.

We introduce a new asymptotic dependence condition.

\begin{definition}
We say that the triangular array $(X_{j,n})_{1 \leq j \leq n, n
\geq 1}$ satisfies {\bf condition (AD-3)} if
$$\lim_{m \to m_1}\limsup_{n \to \infty}n\sum_{j=m+1}^{n}
{\rm Cov}(g(X_{1,n}),g(X_{j,n}))=0, \quad \forall g \in \cC.$$
\end{definition}

\begin{lemma}
For each $n \geq 1$, let $(X_{j,n})_{1 \leq j \leq n}$ be a
strictly stationary sequence of associated random variables with
values in $\BBr \verb2\2 \{0\}$. If the triangular array
$(X_{j,n})_{1 \leq j \leq n, n \geq 1}$ satisfies conditions (AN)
and (AD-3), then it also satisfies condition (AD-1).
\end{lemma}

{\bf Proof:} As in the proof of Lemma \ref{AD1-m-dependent}, it
suffices to show that (\ref{AD1-mdep}) holds. For this, we use the
same ``block'' technique as in the proof of Lemma
\ref{AD1-m-dependent}, except that now the separation blocks have
size $m$ (instead of $m_n$).

For each $1 \leq i \leq k_n$, let $H_{i,n}^{(m)}$ be the (big)
block of consecutive integers between $(i-1)r_n+1$ and $ir_n-m$
and $I_{i,n}^{(m)}$ be the (small) block of size $m$, consisting
of consecutive integers between $ir_n-m+1$ and $ir_n$. Let
$$U_{i,n}^{(m)}=\sum_{j \in H_{i,n}^{(m)}}f(X_{j,n})=N_{ir_n-m,n}(f)-N_{(i-1)r_n,n}(f).$$


Similarly to (\ref{AD-first-conv}) and (\ref{AD-second-conv}), one
can prove that:
\begin{eqnarray}
\label{assoc-first-conv} \lim_{m \to m_1}\limsup_{n \to \infty}|
E(e^{-N_{k_nr_n,n}(f)})-E(e^{-\sum_{i=1}^{k_n}U_{i,n}^{(m)}})| &= &0 \\
\label{assoc-second-conv} \lim_{m \to m_1}\limsup_{n \to \infty}
|\{E(e^{-N_{r_n,n}(f)})\}^{k_n}-\{E(e^{-U_{1,n}^{(m)}})\}^{k_n}|&
=& 0.
\end{eqnarray}

Therefore, in order to prove that relation (\ref{AD1-mdep}) holds,
it suffices to show that
\begin{equation}
\label{AD-third-conv-2} \lim_{m \to m_1} \limsup_{n \to \infty}
|E(e^{-\sum_{i=1}^{k_n}U_{i,n}^{(m)}})-\{E(e^{-U_{1,n}^{(m)}})\}^{k_n}|
=0.
\end{equation}

Without loss of generality, we suppose that the random variables
$(X_{j,n})_{1 \leq j \leq k_n}$ are uniformly bounded. (Otherwise,
we replace them by the random variables $Y_{j,n}=g(X_{j,n}), 1
\leq j \leq k_n$, where $g$ is a bounded non-decreasing function
such that $g(x)=x$ on the support of $f$. The new sequence
$(Y_{j,n})_{1 \leq j \leq k_n}$ consists of uniformly bounded
associated random variables. Moreover, $f(X_{j,n})=f(Y_{j,n})$ for
all $1 \leq j \leq k_n$.)

Moreover, we suppose that the function $f$ satisfies the following
condition: there exists $L_f>0$ such that
\begin{equation}
\label{lipschitz-eq} |f(x)-f(y)| \leq L_f|x-y|, \quad \forall x,y
\in \BBr \verb2\2 \{0\}.
\end{equation} (Note that any function $f
\in C_K^{+}(\BBr \verb2\2 \{0\})$ can be approximated a bounded
sequence of step functions, which in turn can be approximated by a
sequence of functions which satisfy (\ref{lipschitz-eq}).)

Using an induction argument and (\ref{birkel-inequality}), one can
show that:
\begin{equation}
\label{association-induction}
|E(e^{-\sum_{i=1}^{k}U_{i,n}^{(m)}})-\{E(e^{-U_{1,n}^{(m)}})\}^{k}|
\leq L_f^2\sum_{1 \leq i <l \leq k}\sum_{j \in H_{i,n}^{(m)}}
\sum_{j' \in H_{l,n}^{(m)}} {\rm Cov}(X_{j,n},X_{j',n}).
\end{equation}

By stationarity, we have
\begin{eqnarray}
\nonumber \lefteqn{\sum_{1 \leq i <l \leq k_n}\sum_{j \in
H_{i,n}^{(m)}} \sum_{j' \in H_{l,n}^{(m)}} {\rm
Cov}(X_{j,n},X_{j',n}) = \sum_{i=1}^{k_n}(k_n-i) {\rm Cov}(\sum_{j
\in H_{1,n}^{(m)}}
X_{j,n}, \sum_{j' \in H_{i+1,n}^{(m)}}X_{j',n}) =} \\
\label{cov-stationarity} & & \sum_{i=1}^{k_n}(k_n-i) (r_n-m)
\sum_{l=(i-1)r_n+m+1}^{(i+1)r_n-m}{\rm Cov}(X_{1,n}, X_{l,n}) \leq
2 k_nr_n \sum_{l=m+1}^{n}{\rm Cov}(X_{1,n}, X_{l,n}).
\end{eqnarray}

\noindent From (\ref{association-induction}) and
(\ref{cov-stationarity}), we get:
$|E(e^{-\sum_{i=1}^{k_n}U_{i,n}^{(m)}})-\{E(e^{-U_{1,n}^{(m)}})\}^{k_n}|
\leq 2 L_f^2  n \sum_{l=m+1}^{n}{\rm Cov}(X_{1,n}, X_{l,n})$, and
relation (\ref{AD-third-conv-2}) follows from condition (AD-3).
$\Box$

\subsection{Stochastic volatility sequences}

In this subsection, we assume that the dependence structure on
each row of the array $(X_{j,n})_{1 \leq j \leq n, n \geq 1}$ is
that of a {\em stochastic volatility sequence}. More precisely,
\begin{equation}
\label{def-volat-seq} X_{j,n}=\sigma_j Z_{j,n}, \quad 1 \leq j
\leq n, \ n \geq 1,
\end{equation} where $(Z_{j,n})_{1 \leq j \leq n}$ is a
sequence of i.i.d. random variables, and $(\sigma_j)_{j \geq 1}$
is a strictly stationary sequence of positive random variables,
which is independent of the array $(Z_{j,n})_{1 \leq j \leq n, n
\geq 1}$. In this context, $(Z_{j,n})_{1 \leq j \leq n}$ is called
the noise sequence, and $(\sigma_j)_{j \geq 1}$ is called a
volatility sequence. Such models arise in applications to
financial time series (see \cite{davis-mikosch01}). The row
dependence structure among the variables $(X_{j,n})_{1 \leq j \leq
n}$ is inherited from that of the volatility sequence: if
$(\sigma_j)_j$ is $m$-dependent, then so is the sequence
$(X_{j,n})_{1 \leq j \leq n}$. This model is different than a
GARCH model, in which there is a recurrent dependence between the
noise sequence and the volatility sequence.

The dependence structure that we consider for $(\sigma_j)_{j \geq
1}$ is slightly more general than the strongly mixing property.
More precisely, we assume that $(\sigma_j)_{j \geq 1}$ satisfies
the following condition:
\begin{eqnarray*}
(C)  & &  \mbox{there exists a function} \ \psi: \BBn \to
(0,\infty) \ \mbox{with} \ \lim_{m \to \infty} \psi(m)=0, \
\mbox{such that
for any disjoint blocks} \ I,J \\
& & \mbox{of consecutive integers, which are separated by a block
of at least} \ m \ \mbox{integers, and for any} \ (z_j)_j \subset
\BBr
\end{eqnarray*}
\begin{equation}
\label{volatility-cond}  \left|\mbox{Cov}\left(e^{-\sum_{j \in
I}f(\sigma_jz_j)}, e^{-\sum_{j \in J}f(\sigma_jz_j)} \right)
\right| \leq \psi(m), \quad \forall f \in C_{K}^{+}(\BBr \verb2\2
\{0\}).
\end{equation}

We have the following result.

\begin{lemma}
Let $(X_{j,n})_{1 \leq j \leq n, n \geq 1}$ be the triangular
array given by (\ref{def-volat-seq}). If the array $(X_{j,n})_{1
\leq j \leq n, n \geq 1}$ satisfies condition (AN) and
$(\sigma_j)_{j \geq 1}$ satisfies condition (C), then the array
$(X_{j,n})_{1 \leq j \leq n, n \geq 1}$ satisfies condition
(AD-1).
\end{lemma}

{\bf Proof:}
We use the same argument and notation as in the proof of Lemma
\ref{AD1-m-dependent}. Let $(m_n)_n, (r_n)_n$ and $(k_n)_n$ be
sequences of positive integers such that (\ref{properties-of-m-n})
holds, with the function $\psi$ in the place of $\alpha$.

It suffices to prove that (\ref{AD-third-conv}) holds. We now
claim that, for any $1 \leq k \leq k_n$ we have
\begin{equation}
\label{induction-vol}
|E(e^{-\sum_{i=1}^{k}U_{i,n}})-\{E(e^{-U_{1,n}})\}^{k}| \leq (k-1)
\psi(m_n).
\end{equation}
We show this only for $k=2$, the general induction argument being
very similar. Due to the independence between $(\sigma_j)_{j \geq
1}$ and $(Z_{j,n})_{1 \leq j \leq n,  n \geq 1}$, and the
independence of the sequence $(Z_{j,n})_{1 \leq j \leq n}$, we
have:
\begin{eqnarray*}
& & E (e^{-(U_{1,n}+U_{2,n})}) = \int E \left(e^{\textstyle
-\sum_{j\in H_{1,n}}f(\sigma_{j}z_j)-\sum_{j\in
H_{2,n}}f(\sigma_{j}z_j) } \right)
dP(z_j)_{j\in H_{1,n}\cup H_{2,n}} \\
& & = \int \left[E \left(e^{-\sum_{j\in
H_{1,n}}f(\sigma_{j}z_j)-\sum_{i\in H_{2,n}}f(\sigma_{j}z_j)}
\right)-E \left(e^{-\sum_{j\in H_{1,n}}f(\sigma_{j}z_j)}\right) E
\left(e^{-\sum_{j\in
H_{2,n}}f(\sigma_{j}z_j)}\right)\right]dP(z_j)_{
j\in H_{1,n}\cup H_{2,n}} \\
& & + \int E\left(e^{\textstyle -\sum_{j\in
H_{1,n}}f(\sigma_{j}z_j)}\right)dP(z_j)_{j\in  H_{1,n}} \int
E\left(e^{-
\sum_{j\in H_{2,n}}f(\sigma_{j}z_j)}) \right)dP(z_j)_{j\in  H_{2,n}} \\
& & = \int {\rm Cov} \left(e^{-\sum_{j\in
H_{1,n}}f(\sigma_{j}z_j)}, e^{-\sum_{j\in
H_{2,n}}f(\sigma_{j}z_j))}\right) dP(z_j)_{j\in H_{1,n}\cup
H_{2,n}} + E(e^{-U_{1,n}}) E (e^{-U_{2,n}})
\end{eqnarray*}

\noindent where $dP(z_j)_{j\in H_{1,n}\cup H_{2,n}}$ denotes the
law of $(Z_{j,n})_{j\in H_{1,n}\cup H_{2,n}}$, and $dP(z_j)_{j\in
H_{l,n}}$ denotes the law of $(Z_{j,n})_{j\in H_{l,n}}$ for
$l=1,2$. Using condition (\ref{volatility-cond}), we obtain:
\begin{eqnarray*}
|E(e^{-(U_{1,n}+U_{2,n})})-E(e^{-U_{1,n}})E(e^{-U_{2,n}})| & \leq
& \linebreak  \int \left|{\rm Cov}\left(e^{-\sum_{j \in H_{1,n}}
f(\sigma_j z_{j})},e^{-\sum_{j \in H_{2,n}} f(\sigma_j
z_{j})}\right) \right| dP(z_j)_{j\in H_{1,n}\cup H_{2,n}} \\
&\leq &\psi(m_n),
\end{eqnarray*}
since the blocks $H_{1,n}$ and $H_{2,n}$ are separated by a block
of length $m_n$. This concludes the proof of (\ref{induction-vol})
in the case $k=2$.

Relation (\ref{AD-third-conv}) follows using (\ref{induction-vol})
with $k=k_n$, and the fact that $k_n \psi(m_n) \to 0$. $\Box$

\vspace{3mm}

\noindent \small{{\bf Acknowledgement.} The authors would like to
thank an anonymous referee who read the article carefully and made
some suggestions for improving the presentation.}

\normalsize{

\appendix

\section{Poisson Process Representation}

\begin{proposition}
\label{productPiWi} Let $N^{*}$ be a Poisson process on
$(0,\infty)$, with intensity $\rho$. Then $N^*$ admits the
representation $N^* \stackrel{d}{=}\sum_{i \geq
1}\delta_{P_iW_i}$, where $(P_i)_{i \geq 1}$ are the points of a
Poisson process of intensity $\nu$ and $(W_i)_{i \geq 1}$ is an
independent i.i.d. sequence with distribution $F$, if and only if
the measure $\rho$ satisfies:
$$
\int_0^{\infty}(1-e^{-f(x)})\rho(dx)=\int_{0}^{\infty}
\int_{0}^{\infty} (1-e^{-f(wy)})F(dw)\nu(dy)$$ for every
measurable function $f:(0,\infty) \rightarrow (0,\infty)$.
\end{proposition}

{\bf Proof:} Since $N^*$ is a Poisson process of intensity $\rho$,
for any measurable non-negative function $f$, we have
$L_{N^*}(f)=\exp\left\{-\int_{0}^{\infty}(1-e^{-f(x)})\rho(dx)
\right\}$. Let $N^{**}=\sum_{i \geq 1}\delta_{P_iW_i}$. Since $N^*
\stackrel{d}{=}N^{**}$ if and only if $L_{N^{*}}(f)=L_{N^{**}}(f)$
for any measurable non-negative function $f$, the proof will be
complete once we show that
\begin{equation}
\label{Laplace-N**} E\left( e^{-\sum_{i \geq
1}f(P_iW_{i})}\right)=\exp \left\{- \int_{0}^{\infty}
\int_{0}^{\infty}(1-e^{-f(wy)}) F(dw)\nu(dy)\right\}.
\end{equation}

 We first treat the right hand side of (\ref{Laplace-N**}). For this, we
let $g(y)=-\log \int_{0}^{\infty}e^{-f(wy)}F(dw)$. Using the fact
that $F$ is
 a probability measure on $(0,\infty)$ and $M=\sum_{i \geq 1}\delta_{P_i}$ is a Poisson process of
intensity $\nu$, we have
\begin{eqnarray*}
\exp \left\{- \int_{0}^{\infty} \int_{0}^{\infty}(1-e^{-f(wy)})
F(dw)\nu(dy)\right\} &=& \exp \left\{-
\int_{0}^{\infty}(1-e^{-g(y)})\nu(dy)\right\}  =
E\left(e^{-\sum_{i
\geq 1}g(P_i)} \right) \\
&=& E\left(\prod_{i \geq 1}e^{-g(P_i)} \right)
\end{eqnarray*}
For each $i \geq 1$, let $\phi_i(y)=E(e^{-f(W_i y)}), y \geq 0$.
Since $(W_i)_{i \geq 1}$ are i.i.d. random variables with
distribution $F$, for every $y \geq 0$ we have
$\phi_1(y)=\phi_i(y)=\int_{0}^{\infty}e^{-f(wy)}F(dw)=e^{-g(y)},
\forall i \geq 1$.

By considering the random variable $(P_i)_i:\Omega \rightarrow
[0,\infty)^{\BBz_{+}}$ whose law is denoted by $dP(p_i)_{i}$, we
get
\begin{eqnarray*}
\lefteqn{E\left(\prod_{i \geq 1}e^{-g(P_i)} \right)=
E\left(\prod_{i \geq 1}\phi_{i}(P_i) \right)=
\int_{[0,\infty)^{\BBz_{+}}} \prod_{i \geq 1}\phi_{i}(p_i) dP(p_i)_{i} = } \\
& & \int_{[0,\infty)^{\BBz_{+}}} \prod_{i \geq 1}
\left(\int_{\Omega} e^{-f(p_i W_{i}(\omega_{i}))} P(d\omega_i)
\right) dP(p_i)_{i}=\int_{\Omega} \prod_{i \geq 1}
\left(\int_{\Omega} e^{-f(P_i(\omega) W_{i}(\omega_{i}))}
P(d\omega_i) \right) P(d\omega) =\\
& &\int_{\Omega}  \left(\int_{\Omega} \prod_{i \geq 1}
e^{-f(P_i(\omega)W_{i}(\omega'))} P(d\omega') \right) P(d\omega) =
\int_{\Omega}  \prod_{i \geq 1} e^{-f(P_i(\omega)W_{i}(\omega))}
P(d\omega).
\end{eqnarray*}

\noindent For the second last equality above we used the fact that
$(W_{i})_{i \geq 1}$ are independent, whereas for the last
equality above we used the fact that $(P_i)_{i \geq 1}$ and
$(W_{i})_{i \geq 1}$ are independent. $\Box$

\section{Construction of $(r_n)_n$ and $(m_n)_n$ in the
proof of Lemma \ref{AD1-m-dependent}}

\begin{lemma}
If $\lim_{m \to \infty} \alpha(m)=0$, then there exist some
sequences $(r_n)_n$ and $(m_n)_n$ of positive integers such that
$r_n \to \infty$, $k_n:=[n/r_n] \to \infty$, $m_n \to \infty$,
$m_n/r_n \to 0$ and $k_n \alpha(m_n) \to 0$.
\end{lemma}

{\bf Proof:} Denote $\rho_n:=\alpha([\sqrt{n}])$. Clearly $\rho_n
\to 0$. We define
$$\varepsilon_n=\max\{n^{-1/4},\sqrt{\rho_n}\}, \quad \delta_n=\frac{n^{-1/2}}{\varepsilon_n},
 \quad
\eta_n=\frac{\rho_n}{2\varepsilon_n}, \quad
r_n:=[n\varepsilon_{n}], \quad k_n:=[n/r_n], \quad
m_n:=[n\varepsilon_{n} \delta_n]=[\sqrt{n}].$$ 

\noindent Clearly $\varepsilon_n \to 0$ and $m_n \to \infty$. We
will use repeatedly the inequality $x/2 \leq [x] \leq x$, for any
$x \geq 0$. We have: $$r_n \geq \frac{n \varepsilon_n}{2} \geq
\frac{n^{3/4}}{2} \to \infty \quad \mbox{and} \quad k_n \geq
\frac{1}{2} \cdot \frac{n}{r_n} \geq \frac{1}{2} \cdot \frac{n}{n
\varepsilon_n}=\frac{1}{2 \varepsilon_n} \to \infty.$$

\noindent Finally, since $\delta_n \leq n^{-1/4} \to 0$ and
$\eta_n \leq \sqrt{\rho_n}/2 \to 0$, we have:
$$ \quad \frac{m_n}{r_n} \leq \frac{n \varepsilon_n \delta_n}{n
\varepsilon_n/2}=2\delta_n \to 0 \quad \mbox{and} \quad
k_n\alpha(m_n)=k_n\rho_n \leq \frac{n}{r_n} \rho_n \leq \frac{2
}{\varepsilon_n}\rho_n =4\eta_n \to 0.$$
 $\Box$

\end{document}